\documentclass[12pt,a4paper,twoside,final,notitlepage, leqno]{article} 
\usepackage[english]{babel}
\usepackage[T1]{fontenc}  
\usepackage{graphicx}  
\usepackage{amsmath,amsthm,epsfig,amsfonts,bbm}  
\newcommand{\pequationdeb}{$$ \left\{ \begin{minipage}[c]{130mm}}
\newcommand{\pequationfin}{\end{minipage}
                           \right. $$}

\def \smb {{\scriptstyle \bullet }}

\newcommand{\beq}     {\begin{equation}}
\newcommand{\enq}     {\end{equation}}
\newcommand{\be}    {\begin{enumerate}}
\newcommand{\ee}    {\end{enumerate}}

\newcommand{\Bb}

\newcommand{\R}{\mathbb{R}}

\def\tvi {\vrule  height 10pt depth 5pt width 0pt}
\def\tv  {\tvi \vrule}
\def\tvg {\tv ~~}
\def\tvd {~~ \tv}
\def \na{ \noalign {\hrule}  }
\def\hcr {\hfill & \cr}
\def\section*#1{}
\def\resume{\if@twocolumn
\section*{R\'esum\'e}
\else \small
\quotation{\bf \it R\'esum\'e \rule[1mm]{1.5mm}{0.2mm}\vspace{0pt}}
\fi}
\def\endresume{\if@twocolumn\else\endquotation\fi}
\def\abstract{\if@twocolumn
\noindent\section*{{\bf Abstract}}
\else \small
\quotation{\noindent \bf {Abstract.} \rule[1mm]{1.5mm}{0.2mm}\vspace{0pt}}
\fi}
\def\endabstract{\if@twocolumn\else\endquotation\fi}

\begin{document} 
\title{\bf \LARGE      Quartic Parameters \\ ~ \vspace{-1cm}  ~\\ 
for Acoustic Applications   \\ ~ \vspace{-1cm}  ~\\  
 of Lattice Boltzmann  Scheme    ~\\~  }

\author { { \large  Fran\c{c}ois Dubois~$^{ab}$ and   Pierre Lallemand~$^{c}$}     \\ ~\\ 
{\it  \small $^a$   Conservatoire National des Arts et M\'etiers, }  \\
{\it  \small  Department of Mathematics,  Paris, France. }  \\      
{\it \small  $^b$  Department of Mathematics, University  Paris-Sud,} \\
{\it \small B\^at. 425, F-91405 Orsay Cedex, France} \\ 
{\it \small francois.dubois@math.u-psud.fr} \\ 
{\it \small  $^c$  Centre National de la Recherche Scientifique, Paris, France.} \\ %  (retired).} \\  
{\it \small  pierre.lallemand1@free.fr}  \\  ~\\~\\ }  
\date{ 11 June 2011~\protect\footnote{~Invited Presentation, 
Sixth International Conference for Mesoscopic Methods in Engineering and Science
 (ICMMES-2009), Guangzhou City (Canton), Guangdong  Province, China, July 13-17, 2009. 
Published, {\it Computers And Mathematics With Applications}, 
volume 61, pages 3404-3416, june  2011, doi:10.1016/j.camwa.2011.01.011.    }}

\maketitle

\begin{abstract} 
 Using the Taylor expansion method, we show that 
it is possible to improve the lattice Boltzmann method for acoustic applications.  
We derive a formal expansion of the eigenvalues of the discrete approximation
and fit the parameters of the scheme to enforce fourth order accuracy. 
The corresponding discrete equations are solved with the help of symbolic manipulation. 
The solutions are obtained  in the case of D3Q27  lattice Boltzmann scheme. 
Various numerical tests support the coherence of this approach.                    %%%%G
 $ $ \\[5mm]
   {\bf Keywords}: Taylor expansion method,  linearized Navier--Stokes.   
%%%%%%%%%%%%%     $ $ \\[5mm] 
%%%%%%%%%%%%%      {\bf PACS numbers}: 02.70.Ns, 05.20.Dd, 47.10.+g, 47.11.+j.  
\end{abstract}

%%%%%%%%%%%%%%%%%%%%%%%%%%%%%%%%%%%%%%%%%%%%%%%%%%%%%%%%%%%%%%%%%%%%%%%%%%%%%%%%%%%%%
%%%%%%%%%%%%%%%%%%%%%%%%%%%%%%%%%%%%%%%%%%%%%%%%%%%%%%%%%%%%%%%%%%%%%%%%%%%%%%%%%%%%%
\bigskip \bigskip  \newpage \noindent {\bf \large 1) \quad  Introduction}
%%%%%%%%%%%%%%%%%%%%%%%%%%%%%%%%%%%%%%%%%%%%%%%%%%%%%%%%%%%%%%%%%%%%%%%%%%%%%%%%%%%%%
%%%%%%%%%%%%%%%%%%%%%%%%%%%%%%%%%%%%%%%%%%%%%%%%%%%%%%%%%%%%%%%%%%%%%%%%%%%%%%%%%%%%%

\smallskip \noindent $\bullet$ \quad 
The lattice Boltzmann equation methodology is a general framework for approximating 
problems that are modelled by partial differential equations 
under conservative form arising in physics. 
It has been first proposed for fluid dynamics in the context of cellular automata 
(see {\it e.g.} Hardy  {\it et al} \cite{HPP73}, Wolfram  \cite{Wo86}, 
d'Humi\`eres  {\it et al} \cite{DLF86}). The %  now 
classical lattice Boltzmann scheme 
(McNamara-Zanetti  \cite{MZ88},  Higuera {\it et al} \cite{HSB89}, 
Qian  {\it et al} \cite{QDL92}, d'Humi\`eres \cite{DDH92}) has been first developed for
nonlinear fluid problems. It has also the capability to approximate thermal flows.   
(Chen  {\it et al} \cite{COA94}), magnetohydrodynamics (Dellar  \cite{De02}) 
 or coupling between  fluid and  structures (see {\it e.g.} Kwon   \cite{Kw08}).
%  Lallemand and Bouzidi (ref ??)). 

\bigskip  \noindent $\bullet$ \quad 
We have proposed in \cite{Du07, Du08}  the Taylor expansion method
 to analyze formally the d'Humi\`eres lattice Boltzmann scheme  \cite{DDH92} 
when the mesh size tends to zero. We have then replaced the (much more) formal 
Chapman-Enskog expansion methodology used in  \cite{DDH92}  by a simple 
Taylor expansion relative to the grid spacing. % natural numerical geometrical parameter. 
In this way we can obtain  the modified %  equivalent  fd 9 mai 2010
 equations of the scheme            %%%%C
(in the sense of Lerat-Peyret \cite{LP74}  and Warming-Hyett  \cite{WH75};  
see also  \cite{GS86, CBL97, VR99})  at an 
arbitrary order in the general nonlinear case. We made more explicit  
in   \cite{Du09}  the result at third
order accuracy in all  generality. We have observed at this occasion 
that a serious  difficulty with the lattice Boltzmann scheme 
lies in the fact that the equivalent mass conservation equation contains an  %%%%C
{\it a priori}  non-null third order term.
We have also proposed an algorithm to derive the
modified  %  equivalent  fd 9 mai 2010
 equation in the case of  linearized problems  \cite{DL09}. 
We have applied this methodology to derive the so-called ``quartic parameters'' 
to enhance the accuracy of the lattice Boltzmann scheme to simulate   %%%%C   %%%%G
shear waves \cite{DL09}. An application of this result is used   by 
Leriche {\it et al}  \cite{LLL08}  for  computing with  spectral method     %%%%C
and lattice Boltzmann scheme  the eigenmodes of the Stokes problem in a cube. 
This Taylor expansion method can also be used for a detailed analysis of  %%%%C
boundary conditions.
In  \cite{DLT09} we have   shown that the ``magic boundary parameters''
of Ginzburg and Adler \cite{GA94} depend in fact upon the detailed choice of 
the parameters of the 
lattice Boltzmann scheme.  

\bigskip \noindent $\bullet$ \quad 
In this contribution, we consider  linearized athermal acoustics 
in two and three space  %%%%C
dimensions. In Section 2, we recall the essential properties of the d'Humi\`eres   %%%%G
scheme. Then in Section 3, we use the method of formal expansion   %%%%C
 to expand the  discrete
eigenmodes of the acoustic system of partial differential equations 
as the mesh size tends to   %%%%C    %%%%G
zero. In Section 4, we increase the accuracy of the lattice Boltzmann scheme 
with the development of ``quartic'' parameters and develop also a weaker approach by
enforcing isotropy at fourth order accuracy. Applications in two and three space dimensions  %%%%G
(D2Q13 and D3Q27 schemes)  are presented in Sections 5 and 6 respectively.
We detail our version of the D3Q27 scheme and the explicit formulae 
for the determination of the quartic parameters.

%%%%%%%%%%%%%%%%%%%%%%%%%%%%%%%%%%%%%%%%%%%%%%%%%%%%%%%%%%%%%%%%%%%%%%%%%%%%%%%%%%%%%
%%%%%%%%%%%%%%%%%%%%%%%%%%%%%%%%%%%%%%%%%%%%%%%%%%%%%%%%%%%%%%%%%%%%%%%%%%%%%%%%%%%%%
\bigskip \bigskip    \newpage   \noindent {\bf \large 2) \quad  
d'Humi\`eres lattice Boltzmann scheme}
%%%%%%%%%%%%%%%%%%%%%%%%%%%%%%%%%%%%%%%%%%%%%%%%%%%%%%%%%%%%%%%%%%%%%%%%%%%%%%%%%%%%%
%%%%%%%%%%%%%%%%%%%%%%%%%%%%%%%%%%%%%%%%%%%%%%%%%%%%%%%%%%%%%%%%%%%%%%%%%%%%%%%%%%%%%

\noindent $\bullet$ \quad   
In the framework proposed by d'Humi\`eres  \cite{DDH92}, the lattice Boltzmann scheme
uses a regular lattice  $ \, {\cal L } \,$ with typical mesh size $ \, \Delta x \, $  %%%%C
and for each node $ \, x \, $ in  $ \, {\cal L } ,$ 
a discrete set of $\, (J+1) \,$ velocities $ \, {\cal V } \,$ is given. 
In this contribution, the set  $ \, {\cal V } \,$ does not depend on the vertex $x$. 
We introduce a velocity scale $ \, \lambda \,$ and the time step $
\, \Delta t \, $ is obtained according to the so-called ``acoustic scaling'' : 
\begin{equation} \label{acou-saling}  
\Delta t \,=\, {{\Delta x}\over {\lambda}}  \,
\end{equation} 
For  $ \, x \in  {\cal L } \,$ and $\, v_j \in  {\cal V } ,$   the point 
$ \, x + v_j \Delta t \,$ is supposed to be a new vertex of the lattice. 
The unknown of this lattice Boltzmann scheme is the particle distribution 
$ \, f_j(x,\, t) \,$ for   $ \, x \in  {\cal L } ,$ $ \, 0 \leq j \leq J \, $
 and discrete values of time $t$. 
The numerical scheme proceeds into  two major steps.

\smallskip \noindent $\bullet$ \quad 
 The first step describes the relaxation  $  \, f \longrightarrow f^* \,  $ 
of particle distribution $ \, f \,$ towards a locally defined equilibrium. %%%%C
It is local in space and
nonlinear in general.
%%%%but linear for acoustics. 
In this paper devoted to acoustics we consider only linear contributions. %%%%C
It is defined with the help of   a fixed  invertible matrix $ \, M . \, $     %%%%G
 The  moments $ \, m_k \,$ are defined 
through a simple linear relation 
\begin{equation} \label{f-to-m}   
m_k \,=\,   \sum_{j=0}^{J}  M_{k   j} \,  f_j   \,\,, \quad 0 \leq k \leq J \, .  
\end{equation} 
Note that this very simple {\bf linear}  hypothesis %%%%% 7 dec 2010 
is not satisfied in  the  scheme  proposed by Geier  %%%%C
{\it et al}  \cite{GGK06}.
The first $d+1$ moments (where $\, d \,$ is the space dimension, equal to 
$d=2$  or $d=3$  in the present  applications)  are the total density 
$ \, \rho \,$ and the  momentum $ \, q_{\alpha} : \,$ 
\begin{equation} \label{conserves}    \left\{ \begin{array}{rcl} \displaystyle 
m_0  \, =  \, \rho &\, \equiv \,&  
 \displaystyle  \sum_{j=0}^{J} f_j \, \, \\   \displaystyle 
m_{\alpha} \, = \, q_\alpha & \, \equiv \, & 
 \displaystyle  \sum_{j=0}^{J} v_j^{\alpha} \, f_j \,, 
\qquad 1 \leq \alpha \leq d \,.    
\end{array}  \right.    \end{equation}    
We denote by $ \, W \, $ the vector (in $\, \R^N $) composed of the density and the
components of the momentum. These moments are supposed to be at equilibrium: 
\begin{equation} \label{equil-moments}   
m^*_i = m_i\equiv  m^{\rm eq}_i \equiv W_i  \,, \qquad 0 \leq i \leq d \, .
\end{equation} 
The relaxation evolution    $  \, m \longrightarrow m^* \,  $ 
due to linearized collisions is local and trivial for the vector $ W $ of conserved
quantities, as remarked in  (\ref{equil-moments}). For the other moments, a distribution
of equilibrium moments $ \, G(\smb) \,$ is given as a (linear     %%%%C
for the present study of acoustic waves) function of
the vector $W$: 
\begin{equation} \label{equil-gauss}  
 m^{\rm eq}_k \,=\,  G_k(W)  \,, \qquad d <  k \leq J \,. 
\end{equation}  
For $\, k \geq N, \,$ the evolution   $  \, m \longrightarrow m^* \,  $ 
is supposed to be  uncoupled: 
\begin{equation} \label{relax}   
 m^*_k  =  m_k + s_k \,(   m_k^{\rm eq} - m_k) \,, \quad  k > d   . 
\end{equation} 
It is parametrized by the so-called ``relaxation rate parameters''   %%%%C
$\, s_k. $ For stability
reasons (see {\it e.g.} Lallemand and Luo   \cite{LL00}), 
we have the inequalities 
\begin{equation} \label{stabili}   
 0 \,< \, s_k \, < \, 2  \,, \quad  k > d   . 
\end{equation} 
When the components $ \, m^*_k (x, \, t) \,$ are computed for each 
$\, x \in   {\cal L} \,$ at discrete time $\, t,$ the distribution 
$\,  f^*_j (x, \, t) $ after relaxation  is  reconstructed by inversion of  relation 
 (\ref{f-to-m}): 
\begin{equation} \label{m-to-f}  
f_j^* \,=\,   \sum_{\ell=0}^{J}  M^{-1}_{j \ell} \, \, m_{\ell}^*   \,, 
\qquad 0 \leq j \leq J \,.
\end{equation} 

\smallskip \noindent $\bullet$ \quad 
The second step is just a free advective evolution without collision 
through characteristics: 
\begin{equation} \label{advecte}   
 f_j(x,\, t+\Delta t) =  f_j^*(x - v_j \, \Delta t , \, t) \,,  
\qquad 0 \leq j \leq J \,, \quad x \in  {\cal L} \, . 
\end{equation} 

\smallskip \noindent $\bullet$ \quad 
The asymptotic analysis of cellular automata 
provides  evidence supporting   asymptotic   partial differential equations
  with transport coefficients related to the induced          %%%%C
parameter defined by the so-called H\'enon's relation \cite{He87}
\begin{equation} \label{sigma}  
\sigma_k \, \equiv \, {{1}\over{s_k}} \,-\, {1\over2} \, . 
\end{equation} 
The lattice Boltzmann scheme is classically considered as second order accurate 
(see {\it e.g.}  \cite{LL00}).  
We describe in Section 6  the D3Q27 ($d=3$, $J=26$) numerical scheme that we used for 
acoustic application.

%%%%%%%%%%%%%%%%%%%%%%%%%%%%%%%%%%%%%%%%%%%%%%%%%%%%%%%%%%%%%%%%%%%%%%%%%%%%%%%%%%%%%
%%%%%%%%%%%%%%%%%%%%%%%%%%%%%%%%%%%%%%%%%%%%%%%%%%%%%%%%%%%%%%%%%%%%%%%%%%%%%%%%%%%%%
\bigskip \bigskip   \bigskip  \bigskip  \noindent {\bf \large 3) \quad  
 Formal development of eigenmodes }
%%%%%%%%%%%%%%%%%%%%%%%%%%%%%%%%%%%%%%%%%%%%%%%%%%%%%%%%%%%%%%%%%%%%%%%%%%%%%%%%%%%%%
%%%%%%%%%%%%%%%%%%%%%%%%%%%%%%%%%%%%%%%%%%%%%%%%%%%%%%%%%%%%%%%%%%%%%%%%%%%%%%%%%%%%%

\smallskip \noindent $\bullet$ \quad 
With the method described in \cite{DL09}, it is possible to derive the equivalent system
of ($d+1$) equations associated to a lattice Boltzmann scheme.     %%%%C
The system of linearized conservation equations issued from this analysis  can be 
written under the form 
\begin{equation} \label{2.1} 
A(\Delta x, \, \partial) \,\smb\,  W  = 0 \, 
\end{equation} 
which uses  a compact notation  stating that $ \, A(\Delta x, \, \partial)  \,$ 
is a  %%%%C
$\, (d+1) \times (d+1)  \,$ 
($  4 \times 4 \,$ for three-dimensional applications, $ 3 \times 3$  for
two-dimensional models) matrix of  differential operators 
 of high order acting on the conserved    % relative to the conservative 
variables $ W$. 

We seek  the eigenmodes of the operator  $ \,A(\Delta x, \partial) $, 
{\it id est} the eigenvalues $ \, \lambda_j  (\Delta x, \partial) \, $ 
 and the eigenvectors    $ \, r_j  (\Delta x, \partial) \, $ such that 
\begin{equation} \label{2.2} 
A(\Delta x, \, \partial) \, \smb \, r_j  (\Delta x, \partial)  = 
\lambda_j  (\Delta x, \partial)  \,\,  r_j  (\Delta x, \partial)  \,, \quad 0 \leq j \leq d. \, 
\end{equation} 
We introduce the diagonal matrix $ \, \Lambda  (\Delta x, \partial) \,$ composed of
the eigenvalues $\, \lambda_j  (\Delta x, \partial)  \,$ and the square 
matrix  $ \, R  (\Delta x, \partial) \,$ composed of the eigenvectors. Then 
relation  (\ref{2.2}) can be written under the following synthetic  form: 
\begin{equation} \label{2.3} 
A(\Delta x, \, \partial) \, \smb \, R  (\Delta x, \partial)  =  
    R  (\Delta x, \partial) \, \smb \,  \Lambda  (\Delta x, \partial)  \, . \, 
\end{equation} 

\bigskip \noindent $\bullet$ \quad 
Moreover, the operator $ \, A(\Delta x, \partial) \, $ 
in  equation (\ref{2.1})  is a  formal series  %  relatively  to 
in  the variable $ \, \Delta x $
that is truncated at   order 4 in this contribution: 
\begin{equation} \label{2.4} 
A(\Delta x, \, \partial) \equiv 
A_0(\partial)   + \Delta x \, A_1(\partial)  + \Delta x^2 \, A_2(\partial) 
  + \Delta x^3 \, A_3(\partial) + {\rm O}\big( \Delta x^4 \big) \, . 
\end{equation} 
We can apply in  this case the so-called stationary perturbation theory 
 for linear operators 
(see {\it e.g.}  \cite{CDL78} for an elementary introduction). 
First for $\, \Delta x = 0 $, the operator $ \, A_0(\partial) \, $ is 
exactly the perfect linear acoustic model. For $d=3$ (the case $d=2$ is a simple
 consequence of this particular analysis and we omit the details), we have  
\begin{equation} \label{2.41} 
A_0(\partial)  = \begin {pmatrix} \displaystyle    \partial_t  
&  \displaystyle    \partial_x  
&  \displaystyle    \partial_y   
&  \displaystyle    \partial_z  
 \\[3mm]    \displaystyle   c_0^2 \,   \partial_x  
&  \displaystyle    \partial_t   & 0 & 0 
 \\[3mm]    \displaystyle   c_0^2 \,   \partial_y  
&  0 & \displaystyle    \partial_t    & 0 
 \\[3mm]    \displaystyle   c_0^2 \,   \partial_z 
&  0 & 0 & \displaystyle    \partial_t   
\end  {pmatrix} \, . \end{equation} 
Introducing the Laplacian operator $ \,  \displaystyle  
\Delta \equiv   \partial_x^2 +  \partial_y^2  +  \partial_z^2   \,\, $, 
 a system  $ \, R_0  (\partial) \,$   of  reference eigenvectors 
can be given  under the form 
\begin{equation} \label{2.5} 
R_0(\partial)  = \begin {pmatrix} \displaystyle  
           \displaystyle  c_0 \, \sqrt{\Delta}  & 
           \displaystyle   - c_0 \, \sqrt{\Delta}  & 0  & 0  
 \\[3mm]   \displaystyle  c_0^2 \,  \partial_x  &  
           \displaystyle  c_0^2 \,  \partial_x  & 
           \displaystyle   \partial_y  &
           \displaystyle   -\partial_{xz}^2    
 \\[3mm]  \displaystyle   c_0^2 \,  \partial_y  &   
          \displaystyle   c_0^2 \,  \partial_y  &  
          \displaystyle   -\partial_x  &
          \displaystyle   -\partial_{yz}^2     
 \\[3mm]  \displaystyle   c_0^2 \,  \partial_z  &   
          \displaystyle   c_0^2 \,  \partial_z  &   0   &   
     \displaystyle    \partial_x^2  + \partial_y^2    
\end  {pmatrix}  \, .
\end{equation} 
It is elementary to observe the relation  

\begin{equation} \label{2.70} 
A_0(\partial) \, \smb \, R_0  (\partial)  =  
    R_0  (\partial) \, \smb \,  \Lambda_0  (\partial)  \,  
 \end{equation} 
and the associated matrix  $ \,\, \Lambda_0   (\partial) \,\,$ of  order zero is simply 
given by  %%%%C
\begin{equation} \label{2.6} 
\Lambda_0   (\partial)  = \begin {pmatrix} 
 \displaystyle   \partial_t + c_0 \, \sqrt{\Delta}
& 0 & 0 & 0      
\\ 0  &  \displaystyle   \partial_t  - c_0 \, \sqrt{\Delta} & 0 & 0  
 \\  0 & 0 & \displaystyle  \partial_t  & 0 
 \\  0 & 0 & 0 & \displaystyle  \partial_t   
\end  {pmatrix} \, .  \end{equation} 
%%%%%%%%%%  juin 2010  
We observe that even if the eigenvalue $ \, \partial_t \,$ is degenerated, 
{\it i.e.} when the eigenvalues are multiple, the 
family of eigenvectors proposed in  (\ref{2.5}) is complete and defines a basis.
%%%%%%%%%%  fin juin 2010  
This allows us to show the existence of two acoustic waves and two shear waves  %%%%C
(see {\it e.g.}  \cite{LL59}).

\bigskip \noindent $\bullet$ \quad 
The parameter $ \, \Delta x \, $  is supposed to be infinitesimal  
and we introduce a formal expansion of the eigenvalues with {\bf diagonal} matrices 
$ \,  \Lambda_j(\partial) $: 
\begin{equation} \label{2.7} 
\Lambda (\Delta x, \, \partial) \equiv 
\Lambda_0(\partial)   + \Delta x \, \Lambda_1(\partial)  
+ \Delta x^2 \, \Lambda_2(\partial) 
  + \Delta x^3 \, \Lambda_3(\partial) + {\rm O}\big( \Delta x^4 \big) \,  
\end{equation} 
and relative perturbations $ \,  Q_j(\partial) \, $ of the eigenvectors: 
\begin{equation} \label{2.8} 
    R  (\Delta x, \partial)   \equiv R_0(\partial) \,\smb\, \big( {\rm Id} 
  + \Delta x \, Q_1(\partial)  + \Delta x^2 \, Q_2(\partial) 
  + \Delta x^3 \, Q_3(\partial) + {\rm O} ( \Delta x^4 ) \big) \,.   
\end{equation} 
We adopt in this work the 
 scaling condition  for eigenvectors proposed {\it e.g.} by 
Cohen-Tannoudji    {\it   et al} \cite{CDL78}. 
The component of $ \, r_j  (\Delta x, \partial) \,  $ relative to 
the corresponding   unperturbed eigenvector  %  $ \, r_j^0  \,  $ 
is exactly equal to unity. In other words,   all   the diagonal terms 
of the  $ \,  Q_j(\partial) \, $  matrices are null: 
\begin{equation} \label{2.9} 
 Q_j(\partial)_{\, \ell  \ell} = 0 \,, \qquad 0 \leq \ell \leq d  \,, \quad j \geq 1  \,.    
\end{equation} 
We insert the expansions (\ref{2.8}) and  (\ref{2.9}) into the eigenmode condition
 (\ref{2.3}). We introduce the expression $ \, \widetilde{A}_j(\partial) \,$ 
for  the perturbation relative to the unperturbed eigenbasis, {\it id est}    %%%%C
\begin{equation} \label{2.10} 
 \widetilde{A}_j(\partial) \equiv  R_0(\partial)^{-1} \,\smb\, A_j(\partial) 
\,\smb\, R_0(\partial)    \,   .
\end{equation} 
Then we 
identify each term relative to $ \, \Delta x \,$ and we obtain in this manner:  
\begin{equation} \label{2.11}  
 \Lambda_1(\partial)  \equiv   \widetilde{A}_1(\partial) 
+  \Lambda_0(\partial)   \,\smb\,  Q_1(\partial) 
-  Q_1(\partial) \,\smb\, \Lambda_0(\partial)  \end{equation}  
\begin{equation} \label{2.12} 
 \Lambda_2(\partial)   \equiv  \widetilde{A}_2(\partial) 
+  \widetilde{A}_1(\partial)   \,\smb\,  Q_1(\partial) 
-     Q_1(\partial)  \,\smb\,  \Lambda_1(\partial) 
+  \Lambda_0(\partial)   \,\smb\,  Q_2(\partial) 
-     Q_2(\partial)  \,\smb\,  \Lambda_0(\partial) 
\end{equation} 
\begin{equation} \label{2.13}  
 \Lambda_3(\partial)   \equiv    \left\{ \begin{array}{c} \displaystyle 
 \widetilde{A}_3(\partial) 
+  \widetilde{A}_2(\partial)   \,\smb\,  Q_1(\partial) 
-     Q_1(\partial)  \,\smb\,  \Lambda_2(\partial)  \\ 
+  \widetilde{A}_1(\partial)   \,\smb\,  Q_2(\partial) 
-     Q_2(\partial)  \,\smb\,  \Lambda_1(\partial)  \\ 
+  \Lambda_0(\partial)   \,\smb\,  Q_3(\partial) 
-     Q_3(\partial)  \,\smb\,  \Lambda_0(\partial)    
\end{array}  \right.  \, . \end{equation}   
Explicit forms  of the eigenvalues at each order may be found without difficulty   %%%%G
using symbolic manipulation.

\bigskip \noindent $\bullet$ \quad 
We introduce a reference density $\, \rho_0 ,\,$ a reference sound velocity  $\, c_0 ,\,$ 
the shear viscosity $ \, \mu \,$ and the bulk viscosity $ \, \zeta . \, $ 
It is well known (see {\it e.g.} Landau-Lifshitz \cite{LL59}) that the four  
 compressible isothermal  acoustic   linearized  equations take the form:  
\begin{equation} \label{3.1}   \left\{ \begin{array}{rcl} \displaystyle 
\partial_t \rho  + \rho_0 \,  {\rm div} \,{\bf u}  &=&   0  \\[3mm] \displaystyle 
\rho_0 \, \partial_t  u  \,+\,  c_0^2 \, \partial_x \rho    
- \mu \, \Delta u - \Big( \zeta + \frac{\mu}{3}  \Big) \, 
 \partial_x \big(  {\rm div} \,{\bf u} \big)  &=& 0 \\ [3mm] \displaystyle 
\rho_0 \, \partial_t  v   \,+\,  c_0^2 \, \partial_y \rho   
- \mu \, \Delta v -  \Big( \zeta + \frac{\mu}{3}  \Big) \, 
 \partial_y \big(  {\rm div} \,{\bf u} \big)  &=& 0  \\ [3mm] \displaystyle 
\rho_0 \, \partial_t  w   \,+\,  c_0^2 \, \partial_z \rho    
- \mu \, \Delta w -  \Big( \zeta + \frac{\mu}{3}  \Big) \,  
 \partial_z \big(  {\rm div} \,{\bf u} \big)  &=& 0  
\end{array}  \right.    \end{equation}   
with $ \, {\rm div} \, {\bf u} \equiv \partial_x u + \partial_y v + \partial_z w \, $ 
and we recognize  matrix (\ref{2.41}) when $ \,  \mu =   \zeta = 0 .$ 
This system  has  two families of eigenvalues. 

\bigskip \noindent $\bullet$ \quad 
On the one hand, the acoustic waves $ \lambda_\pm \,$ are given by  
\begin{equation} \label{3.2} 
 \lambda_\pm = \partial_t \,- \, \gamma \, \Delta  \, \pm \,  
 c_0   \, \sqrt{\Delta} \, \sqrt{1 + {{\gamma^2}\over{c_0^2}}    \, \Delta  } 
\end{equation}   
and are parametrized by the attenuation of the sound waves $\, \gamma \, $ given according to 
\begin{equation} \label{3.3} 
\gamma ={{1}\over{2 \, \rho_0}} \, \Big( \zeta  +  {{d+1}\over{d}} \, \mu  \Big) \, . 
\end{equation}   
We can expand the previous acoustic wave eigenvalues in powers of $ \gamma $ when the attenuation
is sufficiently small compared to the frequency of the sound waves.  %%%%C    %%%%G
Then we have the expansion 
\begin{equation} \label{3.4} 
 \lambda_\pm =  \partial_t  \, \pm \, c_0  \, \sqrt{\Delta}  - \gamma \, \Delta 
\, \pm \, {{1}\over{2}}  {{\gamma^2}\over{c_0}}   \, \Delta^{3/2}    
%  \, \mp \,  {{1}\over{8}} {{\gamma^4}\over{c_0^3}} \,  \Delta^{5/2}  
%  \,+\, {\rm O}(\gamma^6)    \, . 
 \,+\, {\rm O}(\gamma^4)    \, . 
\end{equation}   

%%%%%%%%%%%%%%%%%%%%%%    juin 2010 
% \bigskip \noindent $\bullet$ \quad 
\noindent 
The previous (relatively unusual)  notations with matrices 
of operators  are {\bf exactly the one used}  when implementing without difficulty 
(due to linearity) our
approach  with a  symbolic manipulation software. 
As is well known, the pseudo-differential operators like $ \,  \sqrt{\Delta} \,$ 
and $ \,  \Delta^{3/2} \,$ are defined by their action on Fourier transforms
(we refer {\it e.g.}  to the book of H\"ormander \cite{Ho85}).
It is also possible to introduce a Fourier decomposition on harmonic waves of the type
$\, {\rm exp} \big( i \, (\omega \, t \,-\, {\bf k}  \smb {\bf x} ) \big) . \,$ 
Then we have the usual change of notation: $ \, \partial_t \equiv i \, \omega ,\,$ 
 $ \, \nabla \equiv -i \,  {\bf k} ,\,$  
$ \, \Delta \equiv -  \mid \!\! {\bf k} \!\! \mid^2 ,\,$ 
$ \, \Delta^{3/2}  \equiv -i   \mid \!\! {\bf k} \!\! \mid^3 ,\,$ {\it etc.} 
The corresponding dispersion relation    associated with  the 
annulation of the eigenvalues $ \lambda_\pm $ given in (\ref{3.4})
allows the introduction of the parameter $ \,\Gamma_{\ell} \,$  
according to 
\begin{equation} \label{ondacous} 
-i \, \omega \, \equiv  \,  \Gamma_\ell  \, =  \, \gamma \, \mid \!\! {\bf k} \!\! \mid^2 
\, \pm \,\,  i \, c_0 \,  \mid \!\! {\bf k} \!\! \mid \, 
\Big( 1 - \frac{\gamma^2}{2\, c_0^2} \mid \!\! {\bf k} \!\! \mid^2 \!  \Big) 
\, + \, {\rm O} \big(  \gamma^4 \, \mid \!\! {\bf k} \!\! \mid^5  \! \big)  \, . 
\end{equation}   
%%%%%%%%%%%%%%%%%%%%%%    fin juin 2010 

\bigskip \noindent $\bullet$ \quad 
On the other hand, the shear wave  $ \, \lambda_0 \,$ define  a 
%%%%%%        degenerate   modif 7 dec 2010
multiple  eigenmode  that is simply given by  %%%%C
\begin{equation} \label{3.5}  
 \lambda_0 =   \partial_t \,- \, \nu \, \Delta \,  
\end{equation}   
%
%%%%%%%%%%%%%%%%%%%%%%    juin 2010 
with $ \, \nu \equiv {{\mu}\over{\rho_0}}. \,$ 
With the spectral notation,  relation (\ref{3.5}) can be written 
with the help of the attenuation $ \,\Gamma_{t} \,$ of the shear waves 
under the form 
\begin{equation} \label{ondecisail} 
- i \, \omega \, \equiv  \,  \Gamma_{t}  \, = \, \nu  \, \mid \!\! {\bf k} \!\! \mid^2  \, .
\end{equation}   
As previously, a family of {\bf two} independent eigenvectors 
with two different ``polarizations'' corresponds  to this
shear eigenvalue.  
Moreover, the perturbation  
$ \,  \Delta x \, A_1(\partial)  + \Delta x^2 \, A_2(\partial) 
 + \Delta x^3 \, A_3(\partial)   + \cdots \,$  
of relation (\ref{2.4}) is {\bf diagonal} 
on the basis of unperturbed eigenvectors proposed at relation  (\ref{2.5}).
Due to well-established   arguments (see {\it e.g.} 
the  textbook of Quantum Mechanics of  
Cohen-Tannoudji    {\it   et al} \cite{CDL78} or {\it e.g.} the sections 
``multiple roots'' and ``degenerate roots'' in the book of Hinch \cite{Hi91}), 
the  expansion   (\ref{2.8})  (\ref{2.9}) is completely valid even 
in this (relatively) complicated degenerate  case.  
%%%%%%%%%%%%%%%%%%%%%%    fin juin 2010 

%%%%%%%%%%%%%%%%%%%%%%%%%%%%%%%%%%%%%%%%%%%%%%%%%%%%%%%%%%%%%%%%%%%%%%%%%%%%%%%%%%%%%
%%%%%%%%%%%%%%%%%%%%%%%%%%%%%%%%%%%%%%%%%%%%%%%%%%%%%%%%%%%%%%%%%%%%%%%%%%%%%%%%%%%%%
\bigskip \bigskip  \bigskip  \bigskip  \noindent {\bf \large 4) \quad  
  Increasing  the accuracy of shear and acoustic waves }
%%%%%%%%%%%%%%%%%%%%%%%%%%%%%%%%%%%%%%%%%%%%%%%%%%%%%%%%%%%%%%%%%%%%%%%%%%%%%%%%%%%%%
%%%%%%%%%%%%%%%%%%%%%%%%%%%%%%%%%%%%%%%%%%%%%%%%%%%%%%%%%%%%%%%%%%%%%%%%%%%%%%%%%%%%%

\smallskip \noindent $\bullet$ \quad 
When we use a lattice Boltzmann scheme, the viscosity $\, \mu \, $ 
admits a representation of the type ({\it c.f.}  the D3Q27 scheme):
\begin{equation} \label{mu-d3q27} 
  \mu \,=\,  {{1}\over{3}} \,   \lambda \,  \Delta x \, \sigma_x  
\end{equation}   
and similarly the bulk viscosity $ \, \zeta \,$ is given for the D3Q27  %%%%C
scheme (to fix the ideas) according to 
\begin{equation} \label{zeta-d3q27} 
  \zeta \,=\,    \lambda \,  \Delta x \, \sigma_e  \, 
\Big( {{5}\over{9}} - c_0^2 \Big) \, 
\end{equation}   
with the notations of Section 6. 
We can identify the developments (\ref{2.7}) on the one hand and 
(\ref{3.4})(\ref{3.5}) on the other hand because, due to 
(\ref{mu-d3q27}), (\ref{zeta-d3q27}) and  (\ref{3.3}), 
the viscosities 
$ \, \gamma \, $ and $ \, \mu \,$ are proportional to $ \, \Delta x .$ 
The usual implementation of lattice Boltzmann scheme consists in adjusting the  %%%%C
eigenvalues $ \, \Lambda_0(\partial)  \,$ and  $ \, \Lambda_1(\partial)  \,$ 
in order to fit the viscosities.

\bigskip \noindent $\bullet$ \quad 
%  For most of the schemes we haves used (D2Q9 \cite{LL00}, D2Q13 \cite{LL03}, 
%  D3Q19  \cite{TKSR02} and D3Q27 \cite{DGKLL}), the fitting of eigenvalues at third order 
%  ({\it id est} the matrix  $ \, \Lambda_2(\partial)  $) do not give any 
%  additional  constraint for the coefficients of the scheme. 
%
When we identify the developments  (\ref{2.7}) and  (\ref{3.5}) 
at the order 4 for the shear waves, we recover the ``quartic parameters''
obtained previously \cite{DL09},  in particular for the D2Q9 \cite{LL00} and 
 D3Q19  \cite{TKSR02} schemes. In these cases, the shear viscosity $ \, \mu \,$ follows
 the ``quartic'' condition  
\begin{equation} \label{q-mu-d2q9} 
  \mu \,=\,  {{1}\over{\sqrt{108}}} \,   \lambda \,  \Delta x   \, . 
\end{equation}   
If we add to the previous study the fitting of the acoustic waves  (\ref{3.4}),  
there is no solution of the equations for the two previous schemes D2Q9 and D3Q19. 
We have then
considered the previous conditions for extended stencils such as the D2Q13 
\cite{Qi93, WB96}   %   \cite{LL03}
and  D3Q27 (see {\it e.g.} \cite{MSYL2k})  %%%%    \cite{DGKLL} 
lattice Boltzmann schemes. 
The results for quartic coefficients are
displayed in Section 6 when fitting all the waves (\ref{3.4}) and (\ref{3.5}) of the 
linear problem. 

\bigskip \noindent $\bullet$ \quad 
In order to take advantage of the flexibility of the d'Humi\`eres scheme, we have
developed an ``isotropic'' condition, not as strong as the quartic condition.  %%%%C
With this isotropic condition,
we do not impose anymore the conditions that the numerical waves   (\ref{2.7}) 
 fit exactly the physical ones  (\ref{3.4}) and (\ref{3.5}) but suggest that 
the numerical waves are isotropic. In other words, the operator matrices 
  $ \, \Lambda_2(\partial)  $ (for the order 3) and 
  $ \, \Lambda_3(\partial)  $ (for the order 4) are differential operators that are
constrained in order to depend only on the radial variable $ \, r \, $ 
($  r^2 \equiv x^2 + y^2 + z^2 $)   and on the operator $ \, \partial / \partial r \, .$ 
When this condition is imposed, it is possible to fit isotropic waves for the four
previous schemes. In Table 1 and 2,  %%%%%%%%%% \ref{t-d2q9-d2q13} and  \ref{t-d3q19-d3q27}, 
we display the
number of ({\it a priori} non-independent) discrete equations that must be satisfied by the
parameters of the lattice Boltzmann schemes when we impose isotropy for the shear waves,
isotropy for the acoustic waves, exact fitting of the shear waves and exact fitting of the
acoustic waves. 
All these discrete equations have been solved exactly with the help of symbolic manipulation. 
To fix the ideas, we give in Section 6 the parametrization of the quartic shear and
acoustic waves for the D3Q27 lattice Boltzmann scheme (solution of 18 discrete equations,
according to Table 2). %%%%%%%%%%  \ref{t-d3q19-d3q27}).

%%%%%%%%%%%%%%%%%%%%%%%%%%%%%%%%%%%%%%%%%%%%%%%%%%%%%%%%%%%%%%%%%%%%%%%%%%%%%%%%%%% 
 
%%%%%%%%%%%%%%%%%%%%%%%%%%%%%%%%%     Table 01      %%%%%%%%%%%%%%%%%%%%%%%%%%%%%% 
\bigskip    \bigskip      \bigskip 
\setbox20=\hbox{ $\,\,$ }
\setbox30=\hbox{ isotropic  shear }
\setbox40=\hbox{ isotropic acoustic }
\setbox50=\hbox{ exact shear }
\setbox60=\hbox{ exact acoustic }
\setbox21=\hbox{ order 3 }
\setbox31=\hbox{ $ 0 $ }
\setbox41=\hbox{ $ 1 $ }
\setbox51=\hbox{ $ 0 $ }
\setbox61=\hbox{ $ 2 $ }
\setbox22=\hbox{ order 4 }
\setbox32=\hbox{ $ 1 $ }
\setbox42=\hbox{ $ 1 $ }
\setbox52=\hbox{ $ 2 $ }
\setbox62=\hbox{ $ 2 $ } 
\setbox66=\hbox{ test  } 
\setbox44=\vbox{\offinterlineskip  \halign {
&\tvg#& # &\tvg#&  # &\tvg#& #  &\tvg#&  #  &\tvg#&  # &\tvd#\cr 
\na&   \box20  &   & \box30 &&  \box40  && \box50  && \box60  \hcr 
\na&   \box21  && \hfill \box31 \hfill  && \hfill \box41 \hfill && \hfill \box51 
\hfill &&  \hfill \box61 \hfill  \hcr 
\na&   \box22  &&  \hfill \box32 \hfill  &&  \hfill  \box42  \hfill  &&  
\hfill \box52  \hfill  &&  \hfill \box62  \hfill  \hcr 
 \na}   }  \centerline{\box44  }
\smallskip \hangindent=7mm \hangafter=1 \noindent {\bf  Table  1}.   \quad {\it 
 Fitting ``isotropic'' and ``exact'' shear and acoustic waves at various orders
  of accuracy. Number of equations for the D2Q9 and D2Q13 lattice Boltzmann schemes.  
Note that the same numbers apply to both  lattices. }
%  \ref{t-d2q9-d2q13}
\bigskip      
%%%%%%%%%%%%%%%%%%%%%%%%%%%%%%%%%%%%%%%%%%%%%%%%%%%%%%%%%%%%%%%%%%%%%%%%%%%%%%%%%%%
 
%%%%%%%%%%%%%%%%%%%%%%%%%%%%%%%%%     Table 02      %%%%%%%%%%%%%%%%%%%%%%%%%%%%%% 
\bigskip          %%%%%%%%     \bigskip     \bigskip  
\setbox20=\hbox{ $\,\,$ }
\setbox30=\hbox{ isotropic  shear }
\setbox40=\hbox{ isotropic acoustic }
\setbox50=\hbox{ exact shear }
\setbox60=\hbox{ exact acoustic } 
\setbox21=\hbox{ order 3 }
\setbox31=\hbox{ $ 0 $ }
\setbox41=\hbox{ $ 1 $ }
\setbox51=\hbox{ $ 0 $ }
\setbox61=\hbox{ $ 2 $ }
\setbox22=\hbox{ order 4 }
\setbox32=\hbox{ $ 10 $ }
\setbox42=\hbox{ $ 2 $ }
\setbox52=\hbox{ $ 13 $ }
\setbox62=\hbox{ $ 3 $ } 
\setbox66=\hbox{ test  } 
\setbox44=\vbox{\offinterlineskip  \halign {
&\tvg#& # &\tvg#&  # &\tvg#& #  &\tvg#&  #  &\tvg#&  # &\tvd#\cr 
\na&   \box20  &   & \box30 &&  \box40  && \box50  && \box60  \hcr 
\na&   \box21  && \hfill \box31 \hfill  && \hfill \box41 \hfill && \hfill \box51 
\hfill &&  \hfill \box61 \hfill  \hcr 
\na&   \box22  &&  \hfill \box32 \hfill  &&  \hfill  \box42  \hfill  &&  
\hfill \box52  \hfill  &&  \hfill \box62  \hfill  \hcr 
 \na}   }  \centerline{\box44  }
\smallskip \hangindent=7mm \hangafter=1 \noindent {\bf  Table  2}.   \quad {\it 
 Fitting ``isotropic'' and ``exact'' shear and acoustic waves at various orders
  of accuracy. Number of equations for the  D3Q19  and D3Q27 lattice Boltzmann schemes. 
Note that the same numbers apply to both  lattices. }
%     \label{t-d3q19-d3q27}   
\bigskip   
 %%%%%%%%%%%%%%%%%%%%%%%%%%%%%%%%%%%%%%%%%%%%%%%%%%%%%%%%%%%%%%%%%%%%%%%%%%%%%%%%%%%

%%%%%%%%%%%%%%%%%%%%%%%%%%%%%%%%%%%%%%%%%%%%%%%%%%%%%%%%%%%%%%%%%%%%%%%%%%%%%%%%%%%%%
%%%%%%%%%%%%%%%%%%%%%%%%%%%%%%%%%%%%%%%%%%%%%%%%%%%%%%%%%%%%%%%%%%%%%%%%%%%%%%%%%%%%%
\bigskip \bigskip  \newpage \noindent {\bf \large 5) \quad  
Linearized acoustics with the D2Q13   scheme}  
%%%%%%%%%%%%%%%%%%%%%%%%%%%%%%%%%%%%%%%%%%%%%%%%%%%%%%%%%%%%%%%%%%%%%%%%%%%%%%%%%%%%%
%%%%%%%%%%%%%%%%%%%%%%%%%%%%%%%%%%%%%%%%%%%%%%%%%%%%%%%%%%%%%%%%%%%%%%%%%%%%%%%%%%%%%

% \bigskip
 \noindent $\bullet$ \quad 
Two numerical ``experiments'' have been performed: emission from the center 
(Figures   1,  2  and 3) 
%% \ref{f-d2q13-diverge}, \ref{f-d2q13-diverge-2}, \ref{f-d2q13-diverge-3})  
and emission from the boundary of the circle (see Figures 1,  4, 5, 6). 
% \ref{f-d2q13-diverge}, \ref{f-d2q13-test2}, \ref{f-d2q13-test1}, 
% \ref{f-d2q13-test3}). 
In both cases, a first order  anti-bounce-back boundary condition  %%%%G
(see {\it e.g.} \cite{bfl}) is implemented 
to impose a given density (pressure) on the boundary:  either a constant when sound is
emitted at the center, or a sinusoidal time dependent  value for the emission by the
boundary. In each case, we measure the value of the density  at a fixed position
in space {\it vs }  time (see Figure 6) %%%%%%%%     \ref{f-d2q13-test3}) 
or we analyze the pressure field at a given time 
 (see Figures   2 to 5). %%%%%%%%  \ref{f-d2q13-diverge-2} to \ref{f-d2q13-test1}). 
Obviously the maximum duration
of the simulation with a radius $R$ and sound velocity $c_0$ is 
$R/c_0$  when the source is at the center and $2R/c_0$ when it is
on the circular boundary. These experiments are sensitive to anisotropic behaviour
of both the velocity and of the attenuation as can be seen in the figures.

%%%%%%%%%%%%%%%%%%%%%%%%%%%%%%%%%     Figure 01      %%%%%%%%%%%%%%%%%%%%%%%%%%%%%% 
\bigskip      \bigskip      \bigskip      
\centerline {\includegraphics[width=.4  \textwidth]{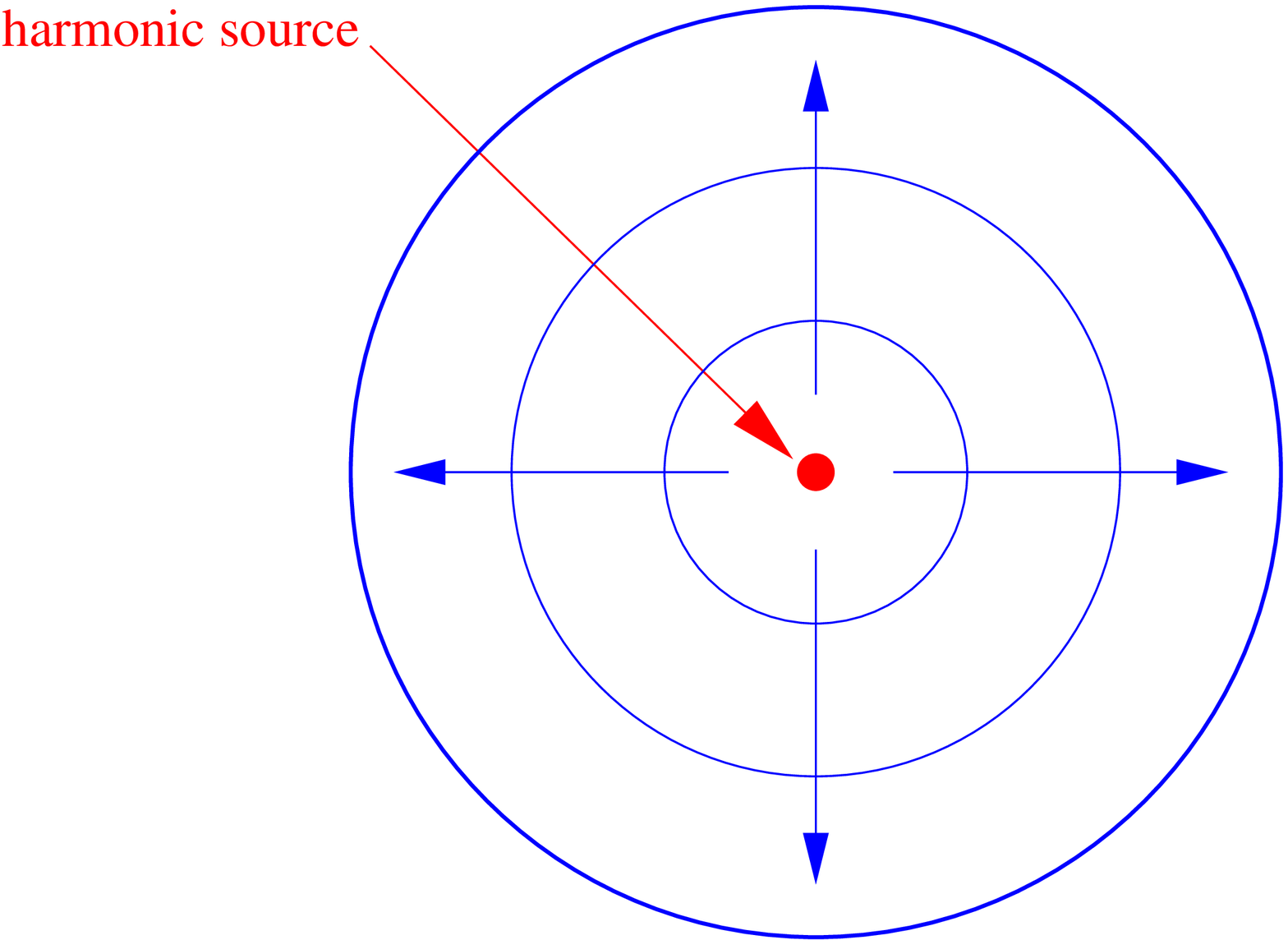} \quad 
 \includegraphics[width=.4  \textwidth]{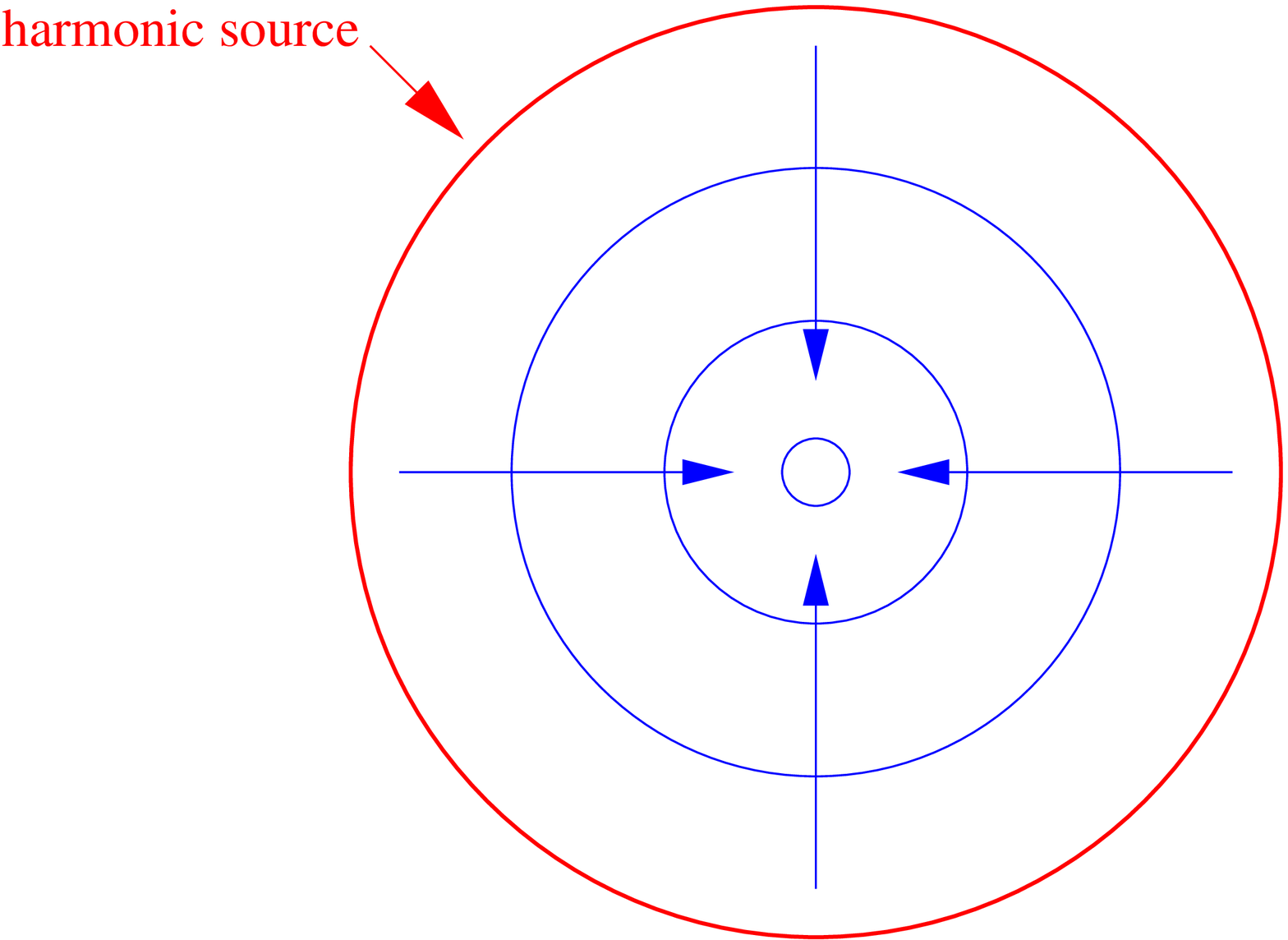}}
\smallskip \hangindent=7mm \hangafter=1 \noindent {\bf Figure 1}.   \quad   {\it  
 Diverging   (left) and converging (right) acoustic test cases.   }
%  \label{f-d2q13-diverge}  \end{figure}   
\bigskip      
%%%%%%%%%%%%%%%%%%%%%%%%%%%%%%%%%%%%%%%%%%%%%%%%%%%%%%%%%%%%%%%%  

%%%%%%%%%%%%%%%%%%%%%%%%%%%%%%%%%     Figure 02      %%%%%%%%%%%%%%%%%%%%%%%%%%%%%% 
\bigskip      \bigskip  \bigskip     
  \centerline { \includegraphics[width=.3  \textwidth]{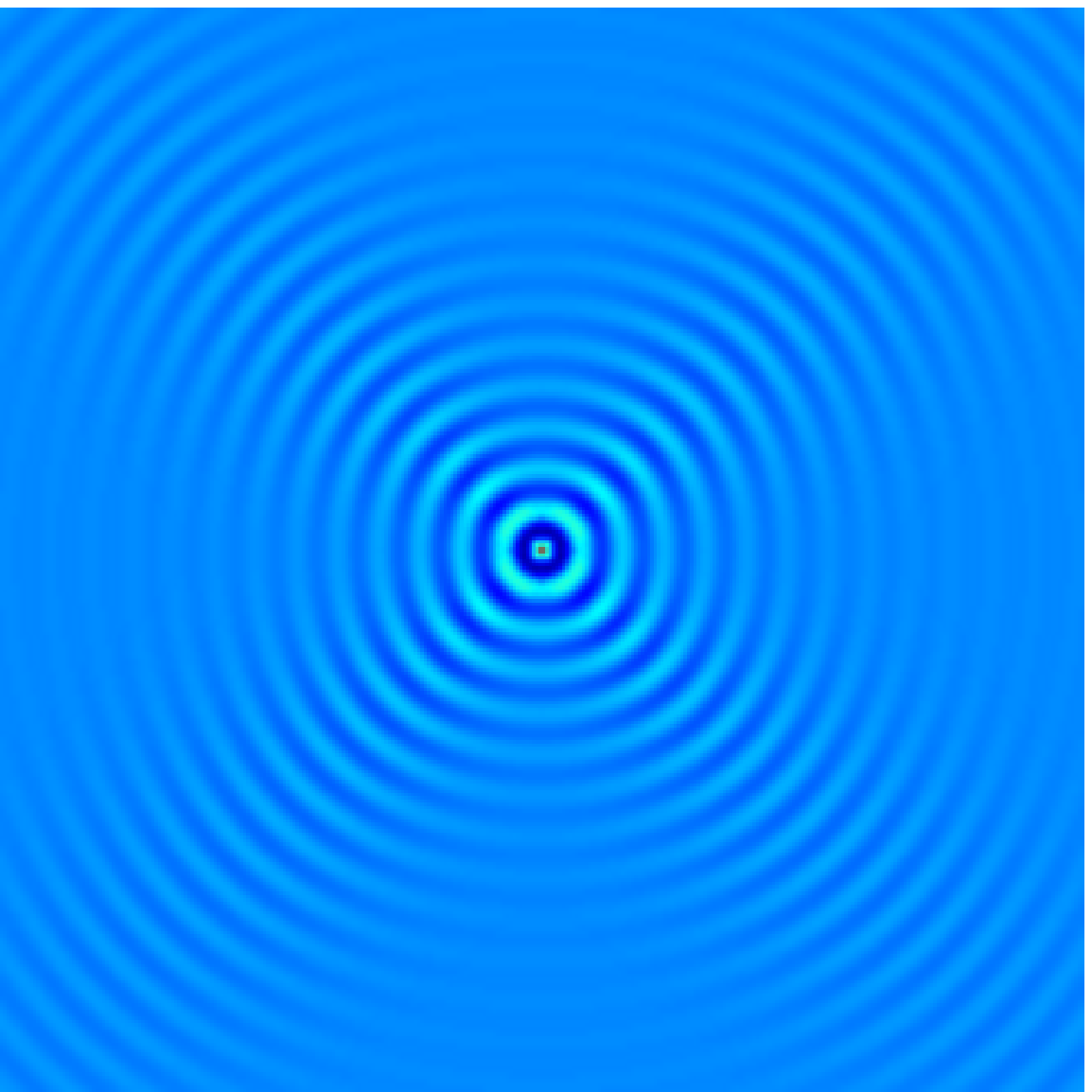}   \quad 
   \includegraphics[width=.3  \textwidth]{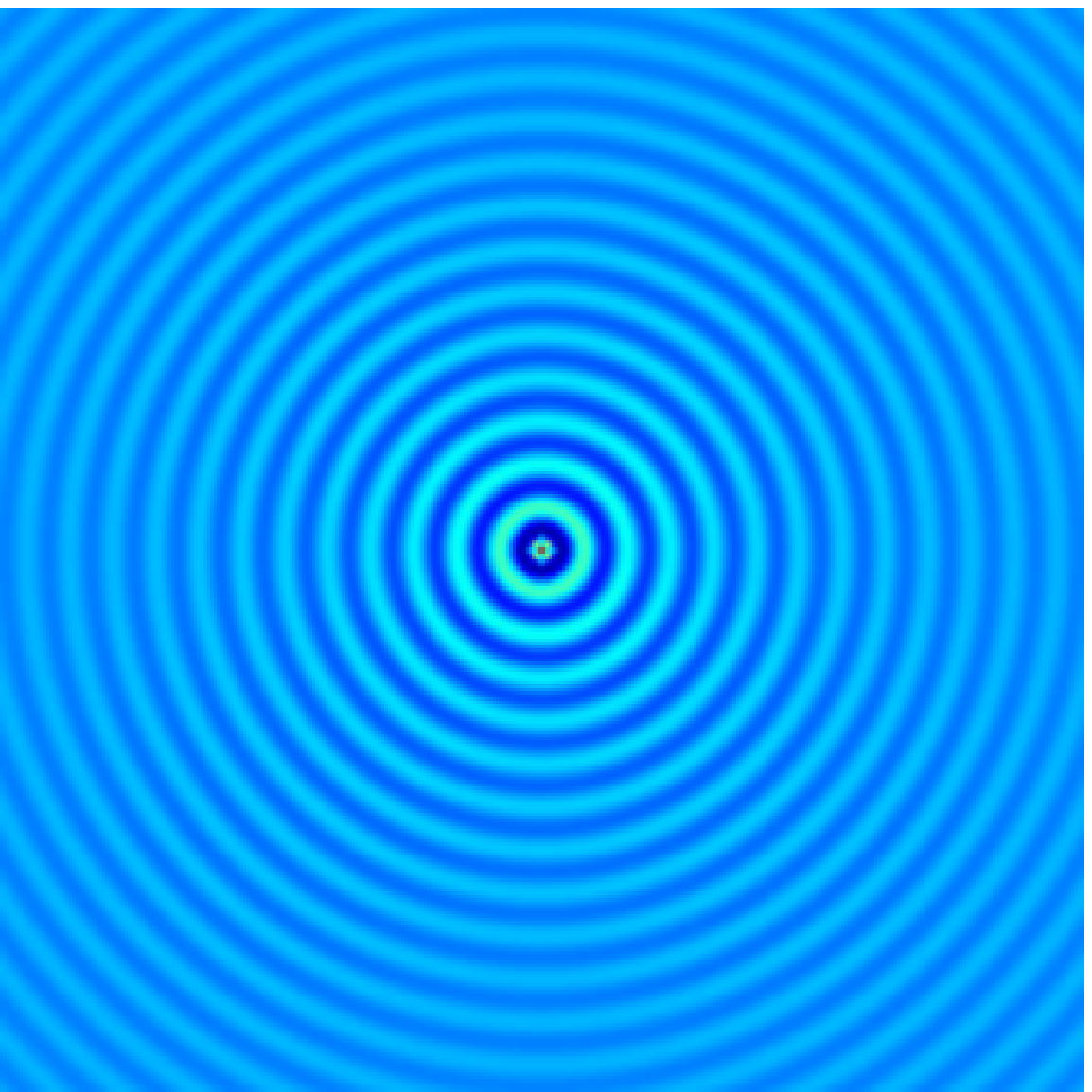} \quad 
   \includegraphics[width=.3  \textwidth]{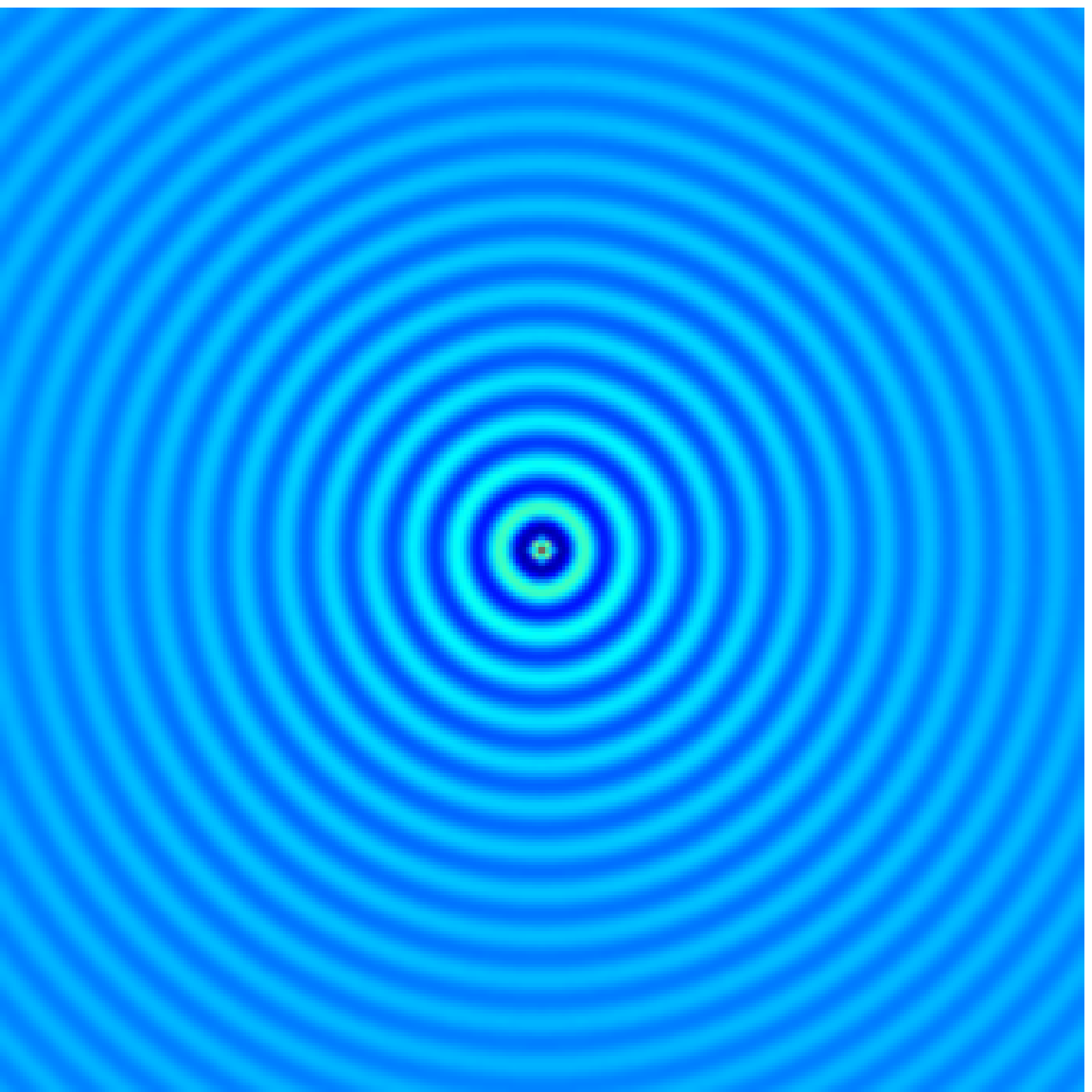}   } 
\smallskip \hangindent=7mm \hangafter=1 \noindent {\bf Figure 2}.   \quad   {\it  
 D2Q13  Acoustic propagation,  diverging acoustic test case, 
time period~=~8 (7~points per wave length, 
 $c_0^2 = 0.80$): usual (left),  isotropic (center) and  quartic   (right).  }
%     \label{f-d2q13-diverge-2}  \end{figure}  
\bigskip      
%%%%%%%%%%%%%%%%%%%%%%%%%%%%%%%%%%%%%%%%%%%%%%%%%%%%%%%%%%%%%%%% 
 
%%%%%%%%%%%%%%%%%%%%%%%%%%%%%%%%%     Figure 03      %%%%%%%%%%%%%%%%%%%%%%%%%%%%%% 
\bigskip     \bigskip  \bigskip      
\centerline { \includegraphics[width=.5  \textwidth]{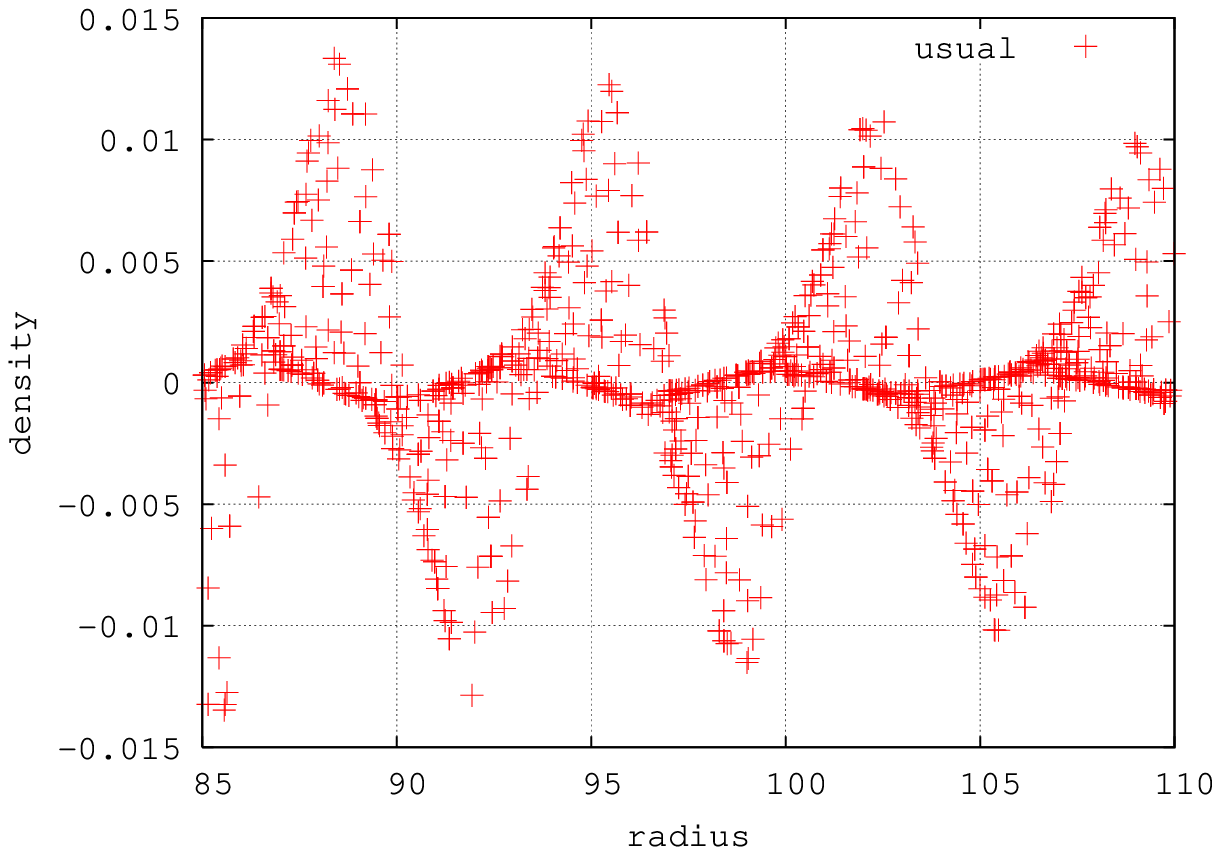} $\,$  
 \includegraphics[width=.5  \textwidth]{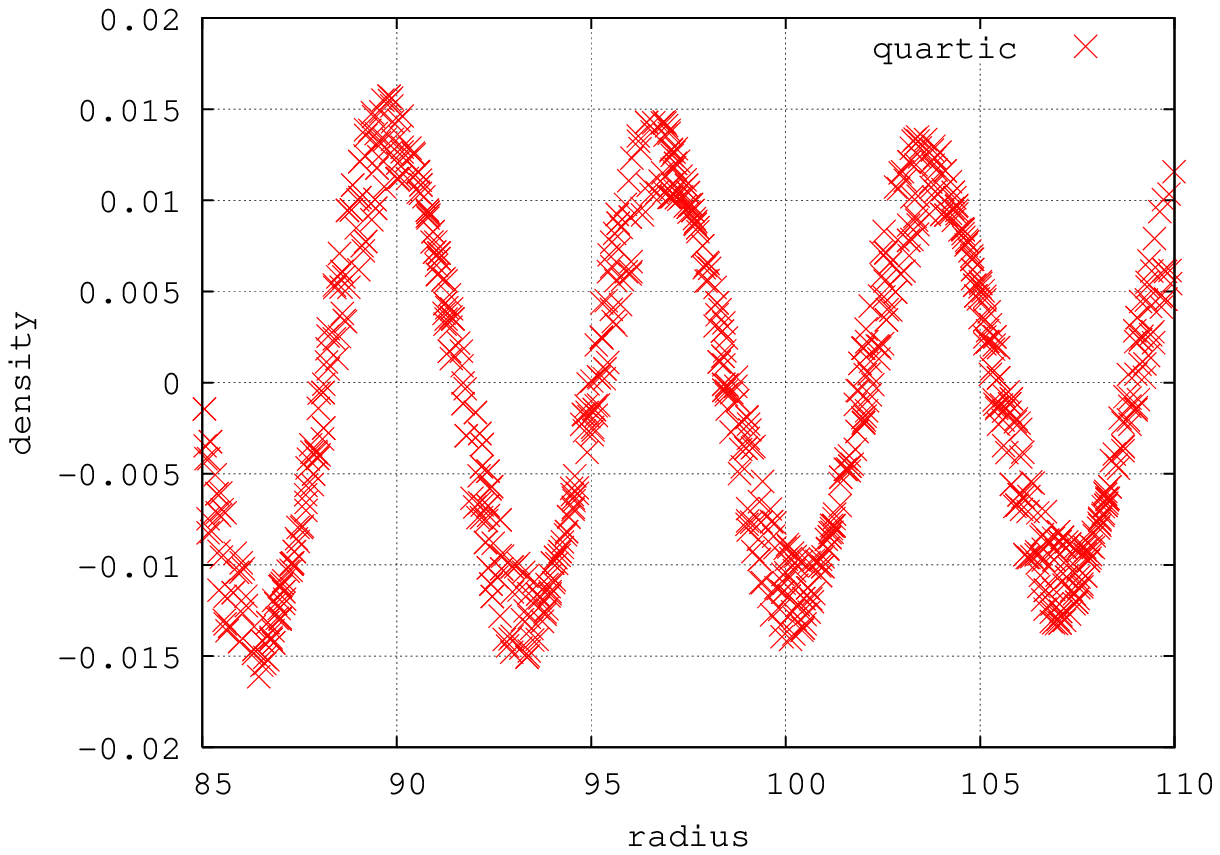}} 
\smallskip \hangindent=7mm \hangafter=1 \noindent {\bf Figure 3}.   \quad   {\it  
 D2Q13  diverging acoustic test case. Measurements after 200 time steps;  
usual parameters (left) and quartic parameters (right).   }  
%     \label{f-d2q13-diverge-3x}  \end{figure}  
\bigskip     
%%%%%%%%%%%%%%%%%%%%%%%%%%%%%%%%%%%%%%%%%%%%%%%%%%%%%%%%%%%%%%%%   

%%%%%%%%%%%%%%%%%%%%%%%%%%%%%%%%%     Figure 04      %%%%%%%%%%%%%%%%%%%%%%%%%%%%%% 
\bigskip      \bigskip  \bigskip      
\centerline { \includegraphics[width=.3  \textwidth]{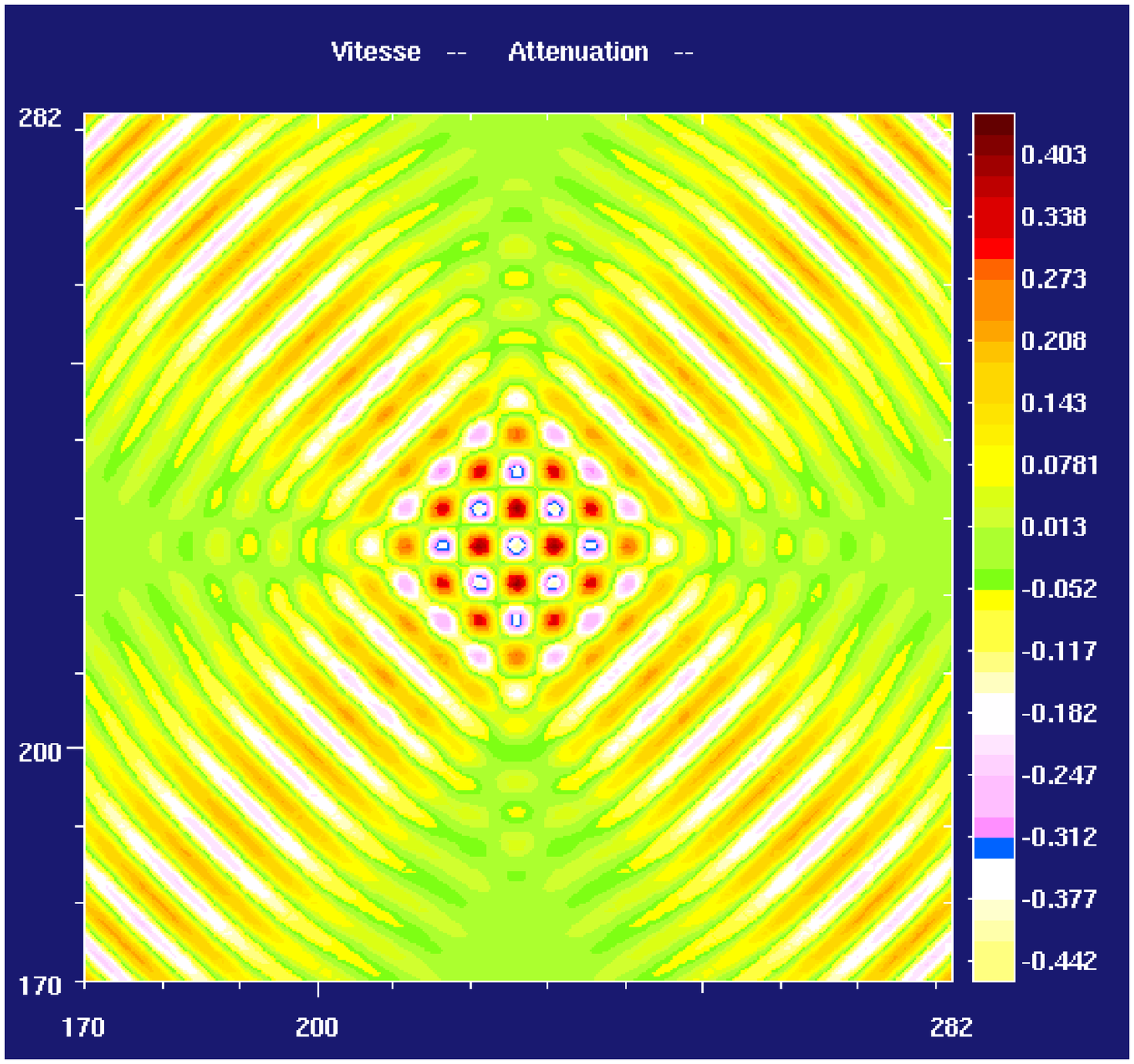}   \quad 
 \includegraphics[width=.3  \textwidth]{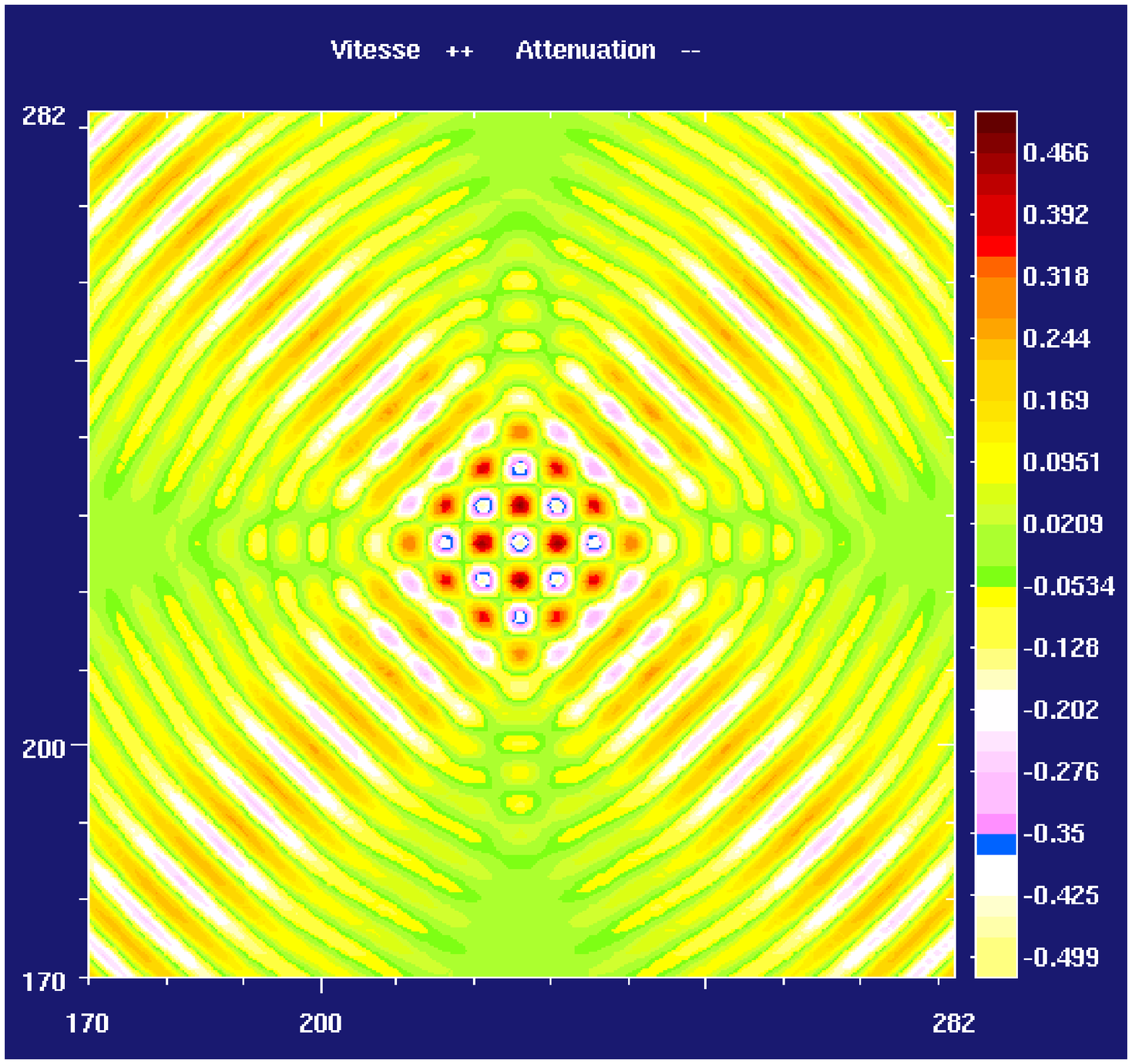} \quad 
 \includegraphics[width=.3  \textwidth]{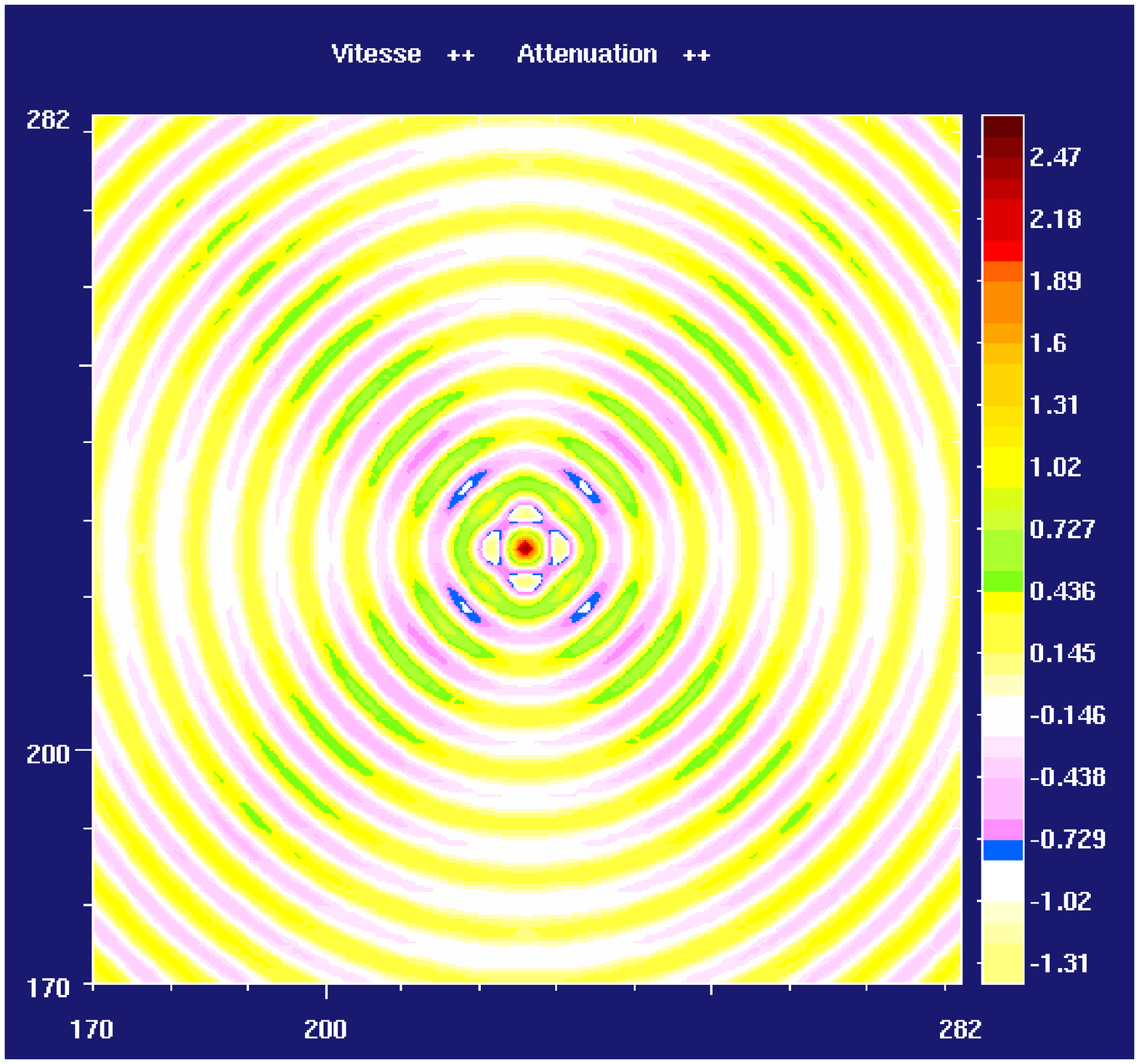}   } 
\smallskip \hangindent=7mm \hangafter=1 \noindent {\bf Figure 4}.   \quad   {\it  
 D2Q13  Converging acoustic propagation, time period = 8 
  (7  points per wave length,  $c_0^2 = 0.80$): 
usual (left),  isotropic (center) and  quartic   (right).  }    
%  \label{f-d2q13-test2}  \end{figure}  
\bigskip     
%%%%%%%%%%%%%%%%%%%%%%%%%%%%%%%%%%%%%%%%%%%%%%%%%%%%%%%%%%%%%%%% 
 
%%%%%%%%%%%%%%%%%%%%%%%%%%%%%%%%%     Figure 05      %%%%%%%%%%%%%%%%%%%%%%%%%%%%%% 
\bigskip \bigskip  \bigskip      
\centerline { \includegraphics[width=.3  \textwidth]{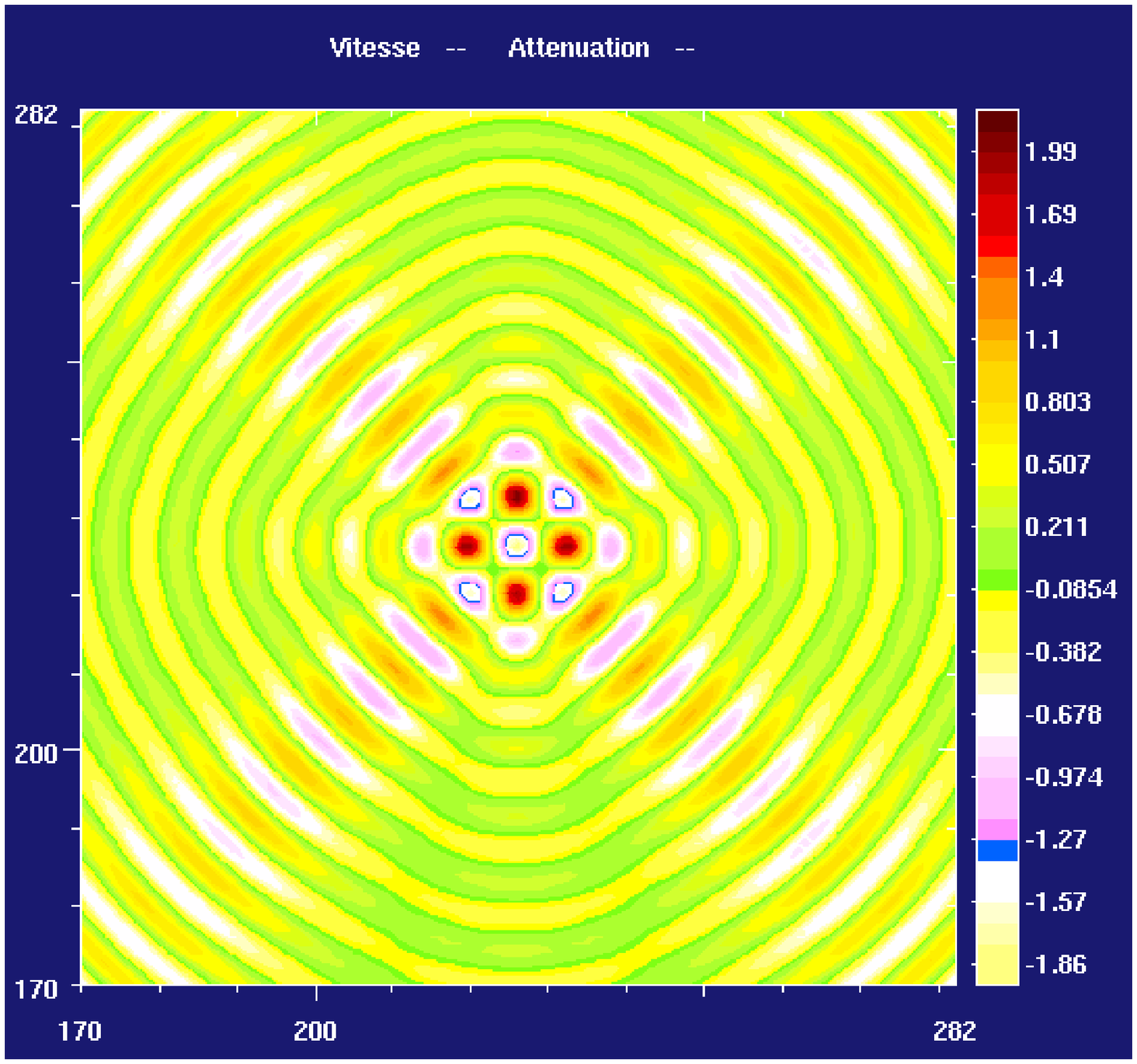}   \quad 
 \includegraphics[width=.3  \textwidth]{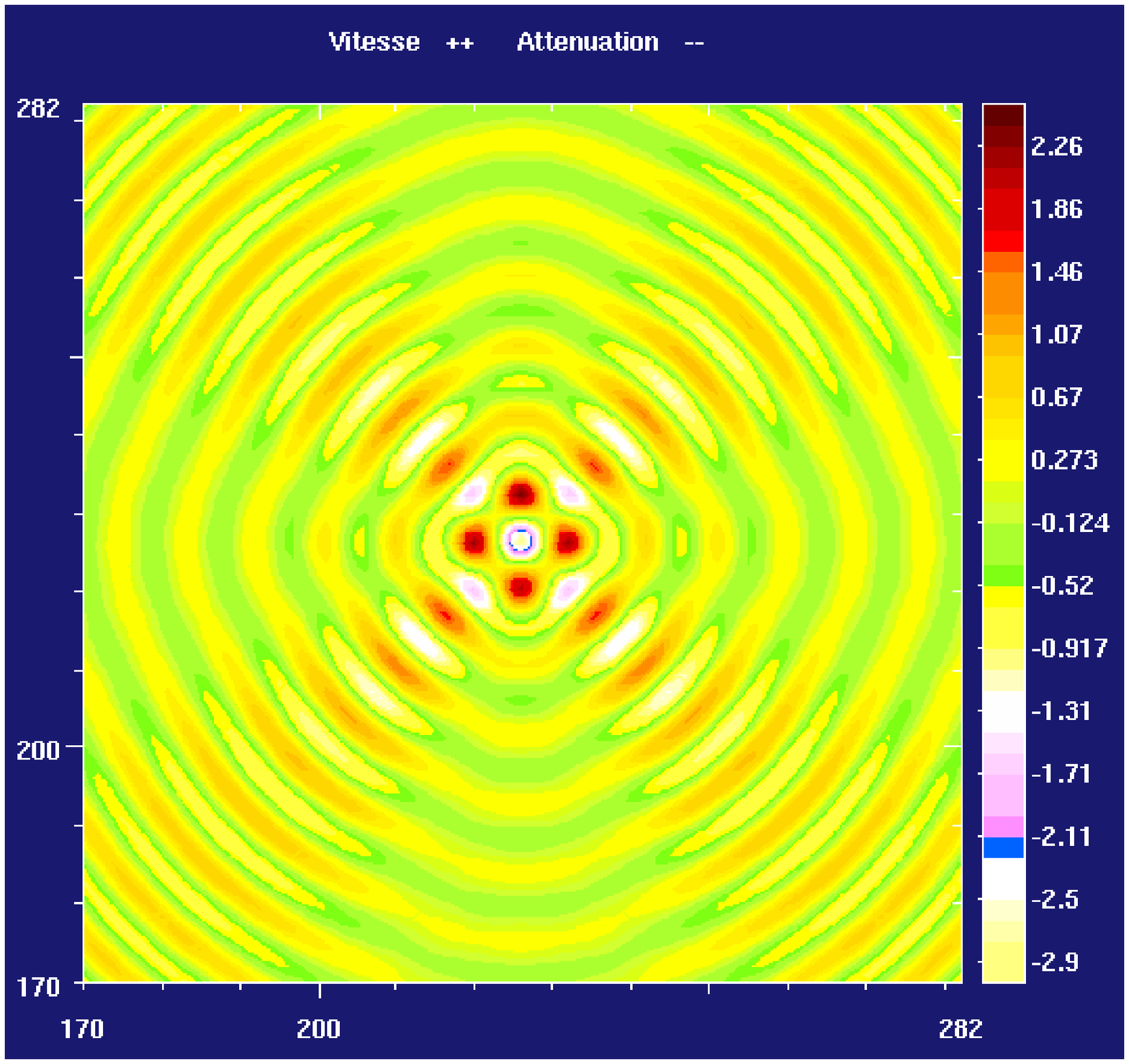} \quad 
 \includegraphics[width=.3  \textwidth]{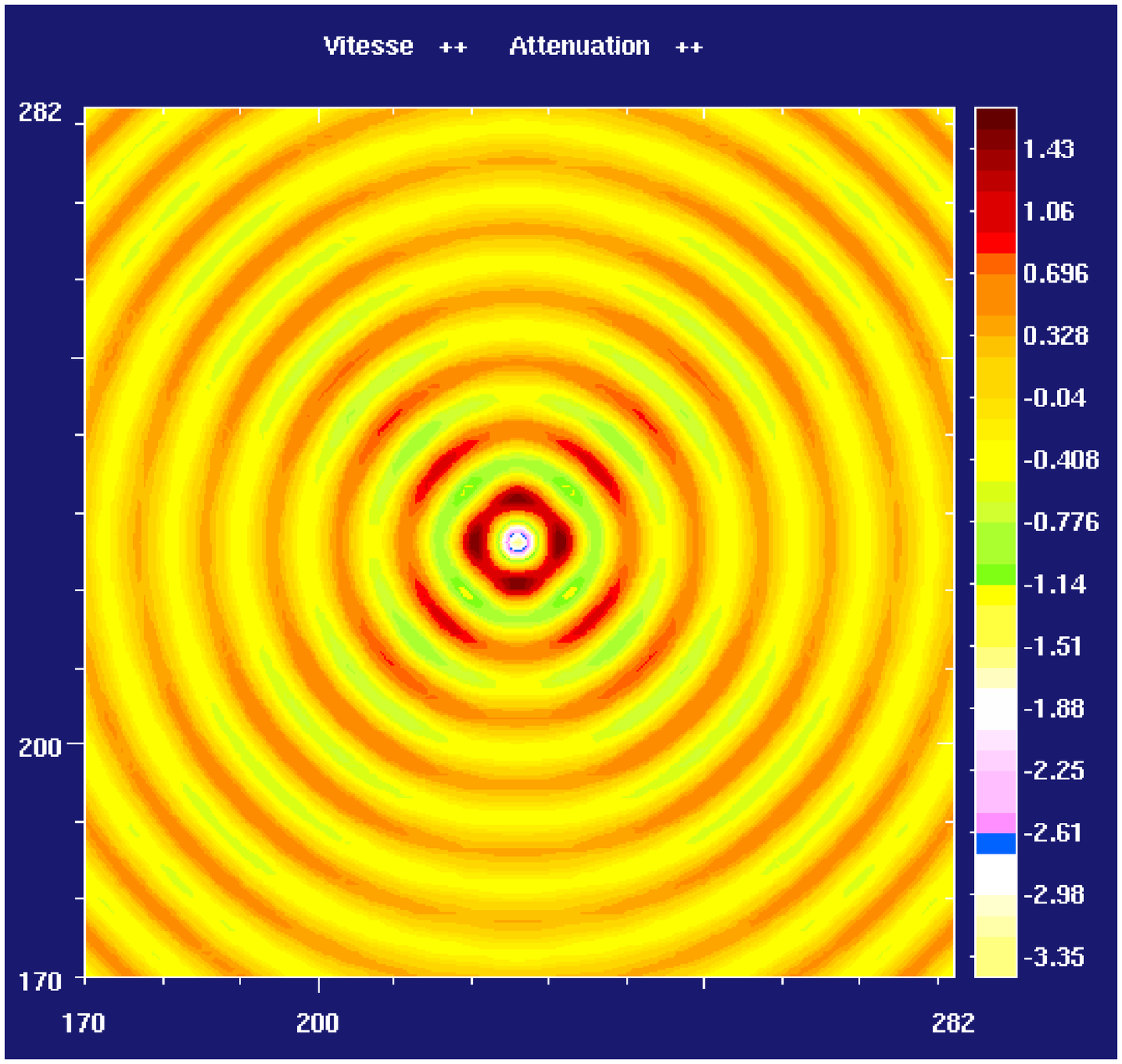}   } 
\smallskip \hangindent=7mm \hangafter=1 \noindent {\bf Figure 5}.   \quad  {\it   
  D2Q13   Converging acoustic propagation, time period = 10  
(9  points per wave  length,  $c_0^2 = 0.80$): 
usual (left),  isotropic (center) and  quartic   (right).   }  
%  \label{f-d2q13-test1}  \end{figure}  
\bigskip     
%%%%%%%%%%%%%%%%%%%%%%%%%%%%%%%%%%%%%%%%%%%%%%%%%%%%%%%%%%%%%%%% 

%%%%%%%%%%%%%%%%%%%%%%%%%%%%%%%%%     Figure 06      %%%%%%%%%%%%%%%%%%%%%%%%%%%%%% 
\bigskip \bigskip  \bigskip      
\centerline { \includegraphics[width=.45  \textwidth]{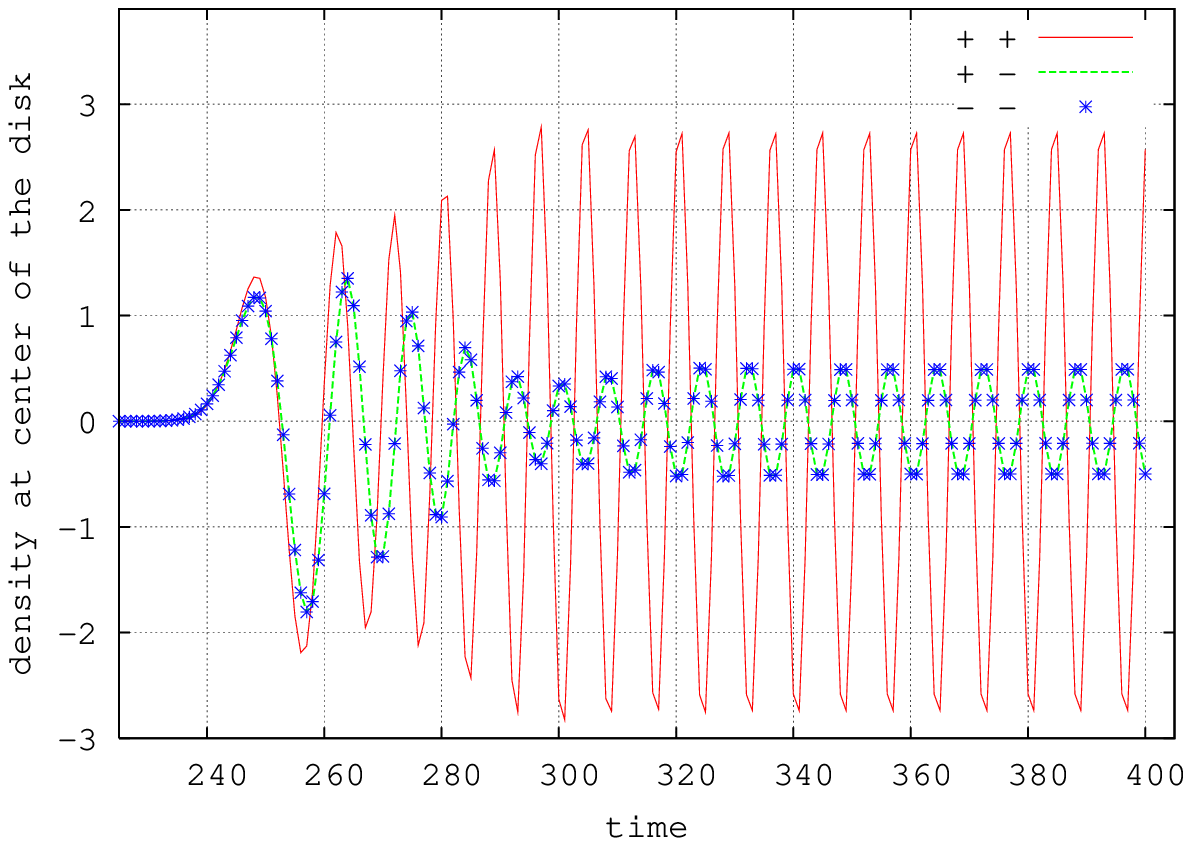}   \quad 
 \includegraphics[width=.45  \textwidth]{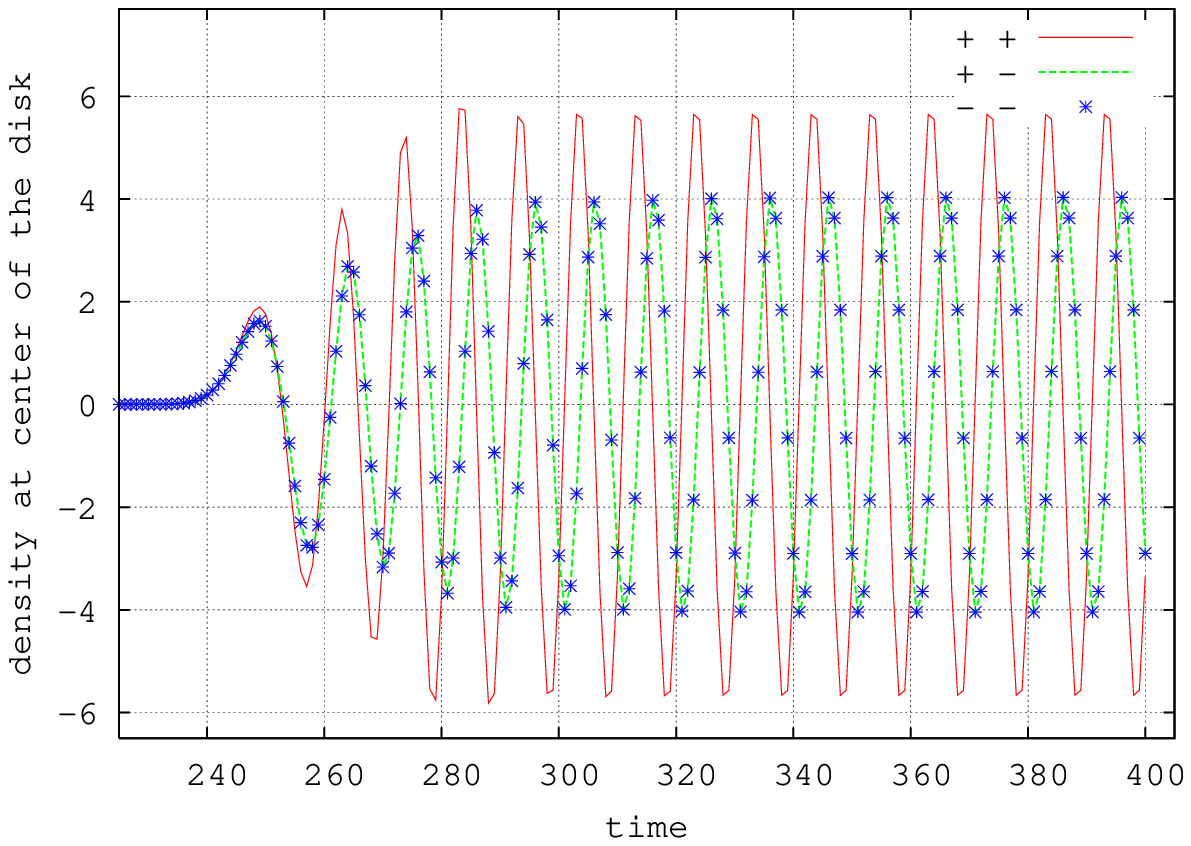}   } 
\smallskip \hangindent=7mm \hangafter=1 \noindent {\bf Figure 6}.   \quad    {\it     
D2Q13   Converging acoustic propagation, time period = 8  (left) and  
time period~=~10   (right).     }
%  \label{f-d2q13-test3}  \end{figure}  
\bigskip     
%%%%%%%%%%%%%%%%%%%%%%%%%%%%%%%%%%%%%%%%%%%%%%%%%%%%%%%%%%%%%%%% 

%%%%%%%%%%%%%%%%%%%%%%%%%%%%%%%%%%%%%%%%%%%%%%%%%%%%%%%%%%%%%%%%%%%%%%%%%%%%%%%%%%%%%
%%%%%%%%%%%%%%%%%%%%%%%%%%%%%%%%%%%%%%%%%%%%%%%%%%%%%%%%%%%%%%%%%%%%%%%%%%%%%%%%%%%%%
\bigskip \bigskip    \noindent {\bf \large 6) \quad D3Q27 athermal   linearized Acoustics}
%%%%%%%%%%%%%%%%%%%%%%%%%%%%%%%%%%%%%%%%%%%%%%%%%%%%%%%%%%%%%%%%%%%%%%%%%%%%%%%%%%%%%
%%%%%%%%%%%%%%%%%%%%%%%%%%%%%%%%%%%%%%%%%%%%%%%%%%%%%%%%%%%%%%%%%%%%%%%%%%%%%%%%%%%%%

\noindent  $\bullet$ \quad  
The D3Q27  Lattice Boltzmann scheme,  described with details 
  {\it e.g.}   in  Mei  {\it et al.}    \cite{MSYL2k} 
%%%%%%%%%       d'Humi\`eres {\it et al.} \cite{DGKLL},   
is simply obtained by taking a tensorial product 
for the set $ \, {\cal V } \,$ of discrete velocities: 
\begin{equation} \label{d3q27-vitesses}  
v_j % ^{\alpha}  
\,=\, \lambda \, \big(  \varepsilon_x \,,\,  \varepsilon_y \,,\,  
\varepsilon_z  \big) \,, \quad 
\varepsilon_x \,,\,   \varepsilon_y ,\, \varepsilon_z \,=\,  -1, 0,  +1  
 \,, \quad  0 \leq j \leq 26 \, . 
\end{equation}  
Due to the large number of moments, we detail in this subsection the way the 
 matrix $M$  parametrizing the transformation (\ref{f-to-m})  is obtained. 
First, velocities $\, v_j^{\alpha} \, $ for 
$\, 0 \leq j \leq J\equiv 26 $ and  $\, 1 \leq \alpha \leq 3 \, $ 
are given according to   relation  (\ref{d3q27-vitesses}).
 The   first four conserved 
moments $\, \rho \, $ and $ \, q^{\alpha}\,$ are determined according to  (\ref{conserves}).
The generation of other moments is done according to the tensorial nature of the  %%%%C
variety of moments that can be constructed, as analyzed by Rubinstein and Luo 
\cite{RL08}: scalar fields are naturally coupled with one another, {\it idem}  for 
vector fields,  and so on. So components of kinetic energy are introduced: 
\begin{equation} \label{en-cinetic} 
\widetilde{m}_4 \, = \, e \, \equiv \, \sum_{j=0}^{26}   
 \Big(   \sum_{\beta=1}^{3} \big( v_j^\beta \big)^2 \Big)  \, f_j  
\end{equation} 
The entire set of second order tensors is completed according to 
\begin{equation} \label{XX-WW}  \left\{ \begin{array}{rccl}  
& \widetilde{m}_5 \, = \,  XX   &  \, \equiv  \, & \displaystyle \sum_{j=0}^{26} 
 \big(  2 \, (v_j^{1})^2 \,-\,  (v_j^{2})^2\,-\,  (v_j^{3})^2 \big)    \, f_j  
\\ 
& \widetilde{m}_6 \, = \,  WW   &  \, \equiv  \, & \displaystyle \sum_{j=0}^{26} 
 \big(  (v_j^{2})^2\,-\,  (v_j^{3})^2 \big)    \, f_j  
\end{array} \right.    \end{equation}  
\begin{equation} \label{XY-YZ-ZX}  \left\{ \begin{array}{rcl}  
& \widetilde{m}_7 \, = \,  XY    \, \equiv  \, & \displaystyle \sum_{j=0}^{26} 
   v_j^{1}  \,  v_j^{2}   \, f_j   
\\ 
& \widetilde{m}_8 \, = \,  YZ    \, \equiv  \, & \displaystyle \sum_{j=0}^{26} 
   v_j^{2}  \,  v_j^{3}   \, f_j    
\\ 
& \widetilde{m}_9 \, = \,  ZX    \, \equiv  \, & \displaystyle \sum_{j=0}^{26} 
   v_j^{3}  \,  v_j^{1}   \, f_j   
\end{array} \right.    \end{equation}  
 Flux of the  energy and  of the square of energy: 
\begin{equation} \label{phi-psi-alpha}    \left\{ \begin{array}{rcl}  
&\widetilde{m}_{9 + \alpha} \,= \, \varphi_{\alpha}  \, \equiv \, &
 \displaystyle  3 \, \sum_{j=0}^{26} \Big(   \sum_{\beta=1}^{3} \big( v_j^\beta \big)^2 
\Big)   \,\, v_j^{\alpha}   \, \, f_j    
\\     
&\widetilde{m}_{12 + \alpha} \, = \,  \psi_{\alpha}     \, \equiv  \, &
 \displaystyle {{9}\over{2}} \, \sum_{j=0}^{26}  
 \Big(   \sum_{\beta=1}^{3} \big( v_j^\beta \big)^2 \Big) ^2 
%  \mid  v_j  \mid ^4  \,   
v_j^{\alpha}   \, \, f_j   
\,,  \qquad 1 \leq \alpha \leq 3 \,    
\end{array} \right.    \end{equation}  
 Square and cube of the energy: 
\begin{equation} \label{e2-e3}    \left\{ \begin{array}{rcl} 
&\widetilde{m}_{16} \, = \, \varepsilon   \, \equiv \, &
 \displaystyle {{3}\over{2}}  \, \sum_{j=0}^{26}  
 \Big(   \sum_{\beta=1}^{3} \big( v_j^\beta \big)^2 \Big) ^2 
%\mid  v_j  \mid ^4  
\,    \, f_j   
\\     
&\widetilde{m}_{17} \, = \,  e_3     \, \equiv  \, &
 \displaystyle {{9}\over{2}}  \, \sum_{j=0}^{26} 
 \Big(   \sum_{\beta=1}^{3} \big( v_j^\beta \big)^2 \Big) ^3  
%  \mid  v_j  \mid ^6  
\,    \, f_j  
\end{array} \right.    \end{equation}  
Product of XX and WW by the energy:
\begin{equation} \label{XXe-WWe}    \left\{ \begin{array}{rcl} 
& \widetilde{m}_{18} \, = \,  XX_e    \, \equiv  \, & \displaystyle 
3 \, \sum_{j=0}^{26}  \big(  2 \, (v_j^{1})^2 \,-\,  (v_j^{2})^2\,-\,  (v_j^{3})^2 \big)
\,  \Big(   \sum_{\beta=1}^{3} \big( v_j^\beta \big)^2 \Big) 
%   \,   \mid  v_j  \mid ^2  
\,  \, f_j  
\\ 
& \widetilde{m}_{19} \, = \,  WW_e    \, \equiv  \, & \displaystyle 
3 \, \sum_{j=0}^{26}   \big(  (v_j^{2})^2\,-\,  (v_j^{3})^2 \big)  \, 
\,  \Big(   \sum_{\beta=1}^{3} \big( v_j^\beta \big)^2 \Big) 
%   \mid  v_j  \mid ^2  
\,  \, f_j    
\end{array} \right.    \end{equation}  
Extra-diagonal second order moments of the energy: 
two different components times energy  and  permutations 
\begin{equation} \label{XYe-YZe-ZXe}  \left\{ \begin{array}{rcl}  
& \widetilde{m}_{20} \, = \,  XY_e    \, \equiv  \, & \displaystyle 
3 \,  \sum_{j=0}^{26}  v_j^{1}  \,  v_j^{2}  \,  %  \mid  v_j  \mid ^2  
\,  \Big(   \sum_{\beta=1}^{3} \big( v_j^\beta \big)^2 \Big) \,\, f_j  
\\
& \widetilde{m}_{21} \, = \,  YZ_e    \, \equiv  \, & \displaystyle 
3 \,  \sum_{j=0}^{26}  v_j^{2}  \,  v_j^{3}  \,   %  \mid  v_j  \mid ^2  
\,  \Big(   \sum_{\beta=1}^{3} \big( v_j^\beta \big)^2 \Big) \,\, f_j  
\\  
& \widetilde{m}_{22} \, = \,  ZX_e    \, \equiv  \, & \displaystyle 
3 \,  \sum_{j=0}^{26}  v_j^{3}  \,  v_j^{1}  \,   %  \mid  v_j  \mid ^2  
\,  \Big(   \sum_{\beta=1}^{3} \big( v_j^\beta \big)^2 \Big) \,\, f_j   
\end{array} \right.    \end{equation}  
Third order pseudovector  $ \,  \tau_{\alpha} $ 
\begin{equation} \label{tau-alpha}  
\widetilde{m}_{22 + \alpha} \, = \,  \tau_{\alpha}    \, \equiv  \,  
\displaystyle    \sum_{j=0}^{26} \, 
v_j^{\alpha} \, \big(  (v_j^{\alpha+1})^2\,-\,  (v_j^{\alpha-1})^2 \big)    \, \,    f_j   \,, 
\qquad 1 \leq \alpha \leq 3 \,.    
\end{equation}  
with a natural permutation convention and the third order totally antisymmetric tensor 
$ \, XYZ $:
\begin{equation} \label{XYZ}  
\widetilde{m}_{26} \, = \,  XYZ    \, \equiv  \,  
\displaystyle    \sum_{j=0}^{26} \, v_j^{1} \, v_j^{2} \,  v_j^{3}  \, \,  f_j \, .   
\end{equation}  
A first (non-orthogonal)  matrix $ \, \widetilde{M} \, $ is obtained by using
 (\ref{conserves})  
 (\ref{en-cinetic}) 
 (\ref{XX-WW}) 
 (\ref{XY-YZ-ZX}) 
 (\ref{phi-psi-alpha}) 
 (\ref{e2-e3}) 
 (\ref{XXe-WWe}) 
 (\ref{XYe-YZe-ZXe}) 
 (\ref{tau-alpha}) and 
 (\ref{XYZ}) 
as linear combinations of the $ \, f_j$'s : 
\begin{equation} \label{M-tilde}  
\widetilde{m}_{k} \, \equiv \,   \sum_{j=0}^{26} \, 
\widetilde{M}_{kj}   \, \,  f_j \,\,, 
\qquad 0 \leq k \leq 26 \,.    
\end{equation}  
Then  matrix $\, M \,$ is orthogonalized from relations
 (\ref{conserves}),
 (\ref{en-cinetic}), 
 (\ref{XX-WW}), 
 (\ref{XY-YZ-ZX}), 
 (\ref{phi-psi-alpha}), 
 (\ref{e2-e3}), 
 (\ref{XXe-WWe}), 
 (\ref{XYe-YZe-ZXe}), 
 (\ref{tau-alpha}) and 
 (\ref{XYZ})    with a Gram-Schmidt orthogonalization   algorithm: 
  \begin{equation*} \label{gram-schmidt}  
 M _{i j}  \,=\,   \widetilde {M} _{i j} - \sum_{\ell < i}  g_{i \ell} 
\, M_{\ell j}  \,,  \quad  i \geq 4    \, .  
 \end{equation*} 
The coefficients $\,  g_{i \ell} \,$ 
 are computed recursively in order to enforce orthogonality: 
\begin{equation*} \label{orthogo}  
\sum_{j = 0}^{26} \, M _{i j} \, M _{k j} \, = \, 0 \, \qquad   \textrm {for } i \not = k  \, .
\end{equation*} 
Note that the coefficients in  
 (\ref{en-cinetic}), 
 (\ref{XX-WW}), 
 (\ref{XY-YZ-ZX}), 
 (\ref{phi-psi-alpha}), 
 (\ref{e2-e3}), 
 (\ref{XXe-WWe}), 
 (\ref{XYe-YZe-ZXe}), 
 (\ref{tau-alpha}) and 
 (\ref{XYZ})  have been chosen such that  % no denominators appear  in matrix $M$. 
the coefficients of matrix $M$ are all integers.  

\smallskip \noindent $\bullet$ \quad 
The equilibrium values of the moments define the equilibrium distribution.
The first moments $ \, \rho \,$ and $ \, q_{\alpha} \, $ are in equilibrium and 
have respectively a scalar and a vectorial structure. 
The equilibrium values of the other moments respect this structure. We have  
\begin{equation} \label{equilibre}    \left\{ \begin{array}{rcl} \displaystyle 
e^{\rm eq} &\,=\,&  \theta \, \lambda^2 \, \rho \\ \displaystyle 
XX^{\rm eq} &\,=\,&  WW^{\rm eq} \,=\,   XY^{\rm eq} \,=\,   YZ^{\rm eq} \,=\, 
  ZX^{\rm eq} \,=\,0   \\ \displaystyle 
\varphi_{\alpha}^{\rm eq} &\,=\,&  c_1 \,  \lambda^2 \, q_{\alpha}   \\ \displaystyle 
\psi_{\alpha}^{\rm eq} &\,=\,&  c_2 \,  \lambda^4 \, q_{\alpha}   \\ \displaystyle 
\varepsilon^{\rm eq} &\,=\,&  \beta  \,  \lambda^4 \, \rho  \\ \displaystyle 
%                                                                     valeurs octobre 2009
%  e_3^{\rm eq} &\,=\,&  \xi  \,  \lambda^4 \, \rho  \\ \displaystyle 
%                                                                 fin valeurs octobre 2009
e_3^{\rm eq} &\,=\,&  \xi  \,  \lambda^6 \, \rho  \\ \displaystyle 
XX_e^{\rm eq} &\,=\,&  WW_e^{\rm eq} \,=\,   XY_e^{\rm eq} \,=\,   YZ_e^{\rm eq} \,=\, 
  ZX_e^{\rm eq} \,=\,0    \\ \displaystyle 
\tau_{\alpha}^{\rm eq} &\,=\,&  c_3 \,  \lambda^2 \, q_{\alpha}   \\ \displaystyle 
XYZ^{\rm eq} &\,=\,& 0 \, . 
\end{array}  \right.    \end{equation}   
We have the following relation between the parameter $ \, \theta \, $ and the 
sound velocity $ \, c_0 : \,$ 
\begin{equation} \label{v-son}  
\theta \,=\,  3 \, c_0^2  - 2   \, . 
\end{equation} 
The relaxation parameters for these moments are denoted respectively by  %%%%C
the following values  
\begin{equation} \label{relaxation}    \left\{ \begin{array}{rcl} \displaystyle 
e^{\rm eq} &\,:\,&  s_e  \\ \displaystyle 
XX^{\rm eq} \,,\,\,   WW^{\rm eq} \,,\,\,  XY^{\rm eq} \,,\,\,  YZ^{\rm eq} \,,\,\,   
  ZX^{\rm eq} &\,:\,&  s_x    \\ \displaystyle 
\varphi_{\alpha}^{\rm eq}  &\,:\,& s_\varphi  \\ \displaystyle 
\psi_{\alpha}^{\rm eq}  &\,:\,& s_\psi   \\ \displaystyle 
\varepsilon^{\rm eq}  &\,:\,& s_\varepsilon  \\ \displaystyle 
e_3^{\rm eq}  &\,:\,&  s_\xi \\ \displaystyle 
XX_e^{\rm eq}  \,,\,\,    WW_e^{\rm eq}    &\,:\,&  s_\gamma   \\ \displaystyle 
XY_e^{\rm eq}  \,,\,\,   YZ_e^{\rm eq}  \,,\,\,   ZX_e^{\rm eq}  
&\,:\,&  s_\chi   \\ \displaystyle 
\tau_{\alpha}^{\rm eq}    &\,:\,&  s_\tau    \\ \displaystyle 
XYZ^{\rm eq}  &\,:\,&  s_\omega  \, . 
\end{array}  \right.    \end{equation}   
With these relaxation parameters, the shear viscosity is isotropic when    %%%%C
\begin{equation} \label{valc1c3}    %%%%C
\left\{     %%%%C
\begin{array}{lcr}    %%%%C
c_1&=&-2 \\    %%%%C
c_3&=&0    %%%%C
\end{array} \right.    %%%%C
\end{equation}    %%%%C

\bigskip \noindent  $\bullet$ \quad   
The quartic parameters for athermal acoustics are solution of the 13 + 2 + 3 = 18 discrete 
equations obtained at Table 2.  %%%%% \ref{t-d3q19-d3q27}. 
They are displayed in this subsection.

\begin{equation} \label{quartic-beta}   \left\{ \begin{array}{rcl}  \displaystyle 
\beta &   \,=\,    \displaystyle     -{{1}\over{2 \,  (\sigma_x-\sigma_e)}}  & \, 
\Big( -9 \, c_0^2 \, \sigma_x
-18 \,  c_0^2  \,  \sigma_e 
+27 \,  c_0^4  \, \sigma_x 
+180 \, c_0^2 \, \sigma_x \, \sigma_e^2 \\ & &  \displaystyle \quad 
+144 \, c_0^2 \, \sigma_x^3  
-8 \, \sigma_x+8 \, \sigma_e
-324  \, c_0^4 \, \sigma_x \, \sigma_e^2  \Big)
\end{array} \right.    \end{equation}    
\begin{equation} \label{quartic-c2}  
c_2  \,=\,  -42 \, \sigma_x^2 + {{5}\over{2}} 
  \end{equation}  
\begin{equation} \label{quartic-sigma-phi}  
\sigma_\varphi   \,=\,  {{1}\over{12 \, \sigma_x}} 
  \end{equation}  
\begin{equation} \label{quartic-sigma-beta}   \left\{ \begin{array}{rcl}  \displaystyle 
\sigma_\varepsilon  & \, = \,  &  \displaystyle   -{{1}\over{48}} \,\, 
{{N_\varepsilon}\over{D_\varepsilon}} \, \\ 
N_\varepsilon & \, = \,  &  \displaystyle 
-76 \, \sigma_x^2
+ 27 \, c_0^2 - 27 \, c_0^4
- 180 \, c_0^2 \, \sigma_e^2
- 468 \, c_0^2 \, \sigma_x^2  
+ 324 \, c_0^4 \, \sigma_e^2
+ 7776  \, c_0^4 \, \sigma_x^4    \\ &&  \displaystyle 
  -93312 \,  c_0^4 \, \sigma_x^4 \, \sigma_e^2
- 4 \, \sigma_e^2 
+ 80 \, \sigma_x \, \sigma_e
- 336 \, \sigma_x^3 \, \sigma_e 
- 1344 \, \sigma_e^2 \, \sigma_x^2
+ 240 \, \sigma_x^4           \\ &&  \displaystyle 
- 10800  \, c_0^2 \, \sigma_x^3 \, \sigma_e 
- 46656 \, c_0^4 \, \sigma_x^3   \, \sigma_e^3
+ 20736 \, \sigma_x^3 \, \sigma_e^3 \, c_0^2 
+ 62208  \, c_0^2 \, \sigma_x^5  \, \sigma_e    \\ &&  \displaystyle 
- 324  \, c_0^4  \, \sigma_x \, \sigma_e         
+ 3888  \, c_0^4 \, \sigma_x \, \sigma_e^3
+ 864 \, c_0^2 \, \sigma_x^2 \, \sigma_e^2 
- 4752 \, c_0^2 \, \sigma_x \, \sigma_e^3      \\ &&  \displaystyle 
+ 51840  \, c_0^2 \, \sigma_e^4  \, \sigma_x^2 
+ 3888 \, c_0^4 \, \sigma_e^2  \, \sigma_x^2    
- 46656 \, c_0^4 \, \sigma_e^4 \, \sigma_x^2 
+ 324  \, c_0^2 \, \sigma_e  \, \sigma_x          \\ &&  \displaystyle 
- 3456 \, \sigma_x^6                     
+ 2880 \, \sigma_x^3 \, \sigma_e^3
+ 20736  \, c_0^2 \, \sigma_x^6 
+ 8064 \, \sigma_x^4 \, \sigma_e^2           
+ 6912 \, \sigma_x^5 \, \sigma_e             \\ &&  \displaystyle 
- 864  \, c_0^2   \, \sigma_x^4              
+ 1440 \, \sigma_x \, \sigma_e^3
- 14400 \, \sigma_e^4 \, \sigma_x^2
+ 324 \, c_0^4 \, \sigma_x^2                  
+ 31104 \, c_0^2 \, \sigma_x^4   \, \sigma_e^2  \\  
D_\varepsilon & \, = \,  &  \displaystyle 
\sigma_x \, (\sigma_x-\sigma_e) \, 
\big( -84 \, \sigma_x^3
+ 432 \, c_0^2\, \sigma_x^3  
+ 84 \, \sigma_x^2 \, \sigma_e
+ 540 \, c_0^2 \, \sigma_x \, \sigma_e^2
- 23 \, \sigma_x                          \\ &&  \displaystyle 
- 972  \, c_0^4  \, \sigma_x \, \sigma_e^2 
- 27 \, c_0^2 \, \sigma_x
+ 81 \, c_0^4 \, \sigma_x  
+ 23 \, \sigma_e
- 54  \, c_0^2 \, \sigma_e  \, \big)
\end{array} \right.    \end{equation}    
We observe that the presence of the factor $ \,  (\sigma_x-\sigma_e) \, $ in the denominator
of the algebraic expression (\ref{quartic-sigma-beta}) 
of $ \, \sigma_\varepsilon \, $ makes this expression incompatible
with a BGK or TRT \cite  {GVD08b}  
 type hypothesis. The same remark holds for   relations 
 (\ref{quartic-sigma-gamma})  and  (\ref{quartic-sigma-chi}) relative to 
the parameters  $ \, \sigma_\gamma \, $ and  $ \, \sigma_\chi \, $ respectively. 
\begin{equation} \label{quartic-sigma-gamma}    \left\{ \begin{array}{rcl}  \displaystyle 
 \sigma_\gamma & \, = \,  &  \displaystyle   {{1}\over{336 \, \sigma_x \, 
(-1+12 \, \sigma_x^2) \, (\sigma_x-\sigma_e)^2}} \,\, N_\gamma \, \\ 
N_\gamma & \, = \,  &  \displaystyle 
1968 \, \sigma_x^4
- 144  \, c_0^2 \, \sigma_x^2 
- 15552  \, c_0^2 \, \sigma_x^4 
- 624 \, \sigma_e^2 \, \sigma_x^2
- 1344 \, \sigma_x^3 \, \sigma_e            \\ &&  \displaystyle 
+ 4608 \, \sigma_x^5 \, \sigma_e                  
+ 8064 \, \sigma_x^4 \, \sigma_e^2
- 12672 \, \sigma_x^6 
+ 103680  \, c_0^2 \, \sigma_x^4  \, \sigma_e^2    \\ &&  \displaystyle 
- 186624 \, c_0^4 \, \sigma_x^4   \, \sigma_e^2
- 27 \, c_0^4                                     
+ 27 \, c_0^2                                    
- 4 \, \sigma_e^2
+ 324 \, c_0^4 \, \sigma_e^2
- 180 \, c_0^2 \, \sigma_e^2                    \\ &&  \displaystyle 
+ 80 \, \sigma_x \, \sigma_e
- 76 \, \sigma_x^2
+ 82944  \, c_0^2   \, \sigma_x^6                 
+ 15552  \, c_0^4  \,\sigma_x^4 
\end{array} \right.    \end{equation}  
\begin{equation} \label{quartic-sigma-chi}     \left\{ \begin{array}{rcl}  \displaystyle 
 \sigma_\chi & \, = \,  &  \displaystyle   {{1}\over{96 \, \sigma_x \, 
(-1 + 84  \, \sigma_x^2) \, (\sigma_x-\sigma_e)^2}} \,\, N_\chi \, \\ 
N_\chi & \, = \,  &  \displaystyle
192 \, \sigma_x^4
- 828 \, c_0^2 \, \sigma_x^2  
+ 972 \, c_0^4 \, \sigma_x^2
+ 6480 \, c_0^2 \, \sigma_x^2 \, \sigma_e^2
- 2592  \, c_0^2  \, \sigma_x^4                   \\ &&  \displaystyle 
- 11664 \, c_0^4 \, \sigma_e^2  \, \sigma_x^2
+ 192 \, \sigma_e^2 \, \sigma_x^2 
- 384 \, \sigma_x^3 \, \sigma_e
+ 2304 \, \sigma_x^5 \, \sigma_e 
+ 4032 \, \sigma_x^4 \, \sigma_e^2            \\ &&  \displaystyle 
- 6336 \, \sigma_x^6
+ 76 \, \sigma_x^2
+ 51840  \, c_0^2 \, \sigma_x^4  \, \sigma_e^2
- 93312 \, c_0^4 \, \sigma_x^4   \, \sigma_e^2  
+ 27 \, c_0^4
- 27 \, c_0^2                                    \\ &&  \displaystyle 
%%               coquille reperee ligne suivante par Pierre le 17 mai 2010 
%%  + 4  \, c_0^4 \, \sigma_e^2 - 324 \, \sigma_e^2     
+ 4   \, \sigma_e^2 - 324 \,   c_0^4 \, \sigma_e^2     %%   coquille reperee 
+ 180 \, c_0^2 \, \sigma_e^2
- 80 \, \sigma_x \, \sigma_e
+ 41472 \, c_0^2 \, \sigma_x^6  
+ 7776  \, c_0^4  \, \sigma_x^4 
\end{array} \right.    \end{equation}   
\begin{equation} \label{quartic-sigma-tau}     \left\{ \begin{array}{rcl}  \displaystyle  
\sigma_\tau & \, = \,  &  \displaystyle   {{1}\over{12}} \,\, {{N_\tau}\over{D_\tau}} \, \\ 
N_\tau & \, = \,  &  \displaystyle 
76 \, \sigma_x^2 
- 324 \, c_0^4 \, \sigma_e^2
+ 180 \, c_0^2 \, \sigma_e^2
- 144 \, \sigma_x^4
- 80 \, \sigma_x \, \sigma_e    
+ 3456  \, c_0^2 \, \sigma_x^4                \\ &&  \displaystyle 
+ 720 \, \sigma_e^2 \, \sigma_x^2
- 576 \, \sigma_x^3 \, \sigma_e 
- 504  \, c_0^2  \, \sigma_x^2 
+ 4320 \, c_0^2 \, \sigma_x^2 \, \sigma_e^2    
+ 27 \, c_0^4
+ 4 \, \sigma_e^2 \\ &&  \displaystyle 
- 27 \, c_0^2
- 7776  \, c_0^4 \, \sigma_e^2  \, \sigma_x^2
+ 648 \, c_0^4 \, \sigma_x^2 \\  
D_\tau & \, = \,  &  \displaystyle  \sigma_x \, \big( 
76 \, \sigma_x^2 
- 528 \, \sigma_x^4
- 504  \, c_0^2 \, \sigma_x^2 
+ 648 \, c_0^4 \, \sigma_x^2
+ 4320 \, c_0^2 \, \sigma_x^2 \, \sigma_e^2    \\ &&  \displaystyle 
+ 3456 \, c_0^2 \, \sigma_x^4  
- 7776  \, c_0^4  \, \sigma_e^2 \, \sigma_x^2
+ 336 \, \sigma_e^2 \, \sigma_x^2
+ 192 \, \sigma_x^3 \, \sigma_e
+ 27 \, c_0^4
- 27 \, c_0^2    \\ &&  \displaystyle 
+ 4 \, \sigma_e^2 - 324 \, c_0^4 \, \sigma_e^2
+ 180 \, c_0^2 \, \sigma_e^2
- 80 \, \sigma_x \, \sigma_e \big) 
\end{array} \right.    \end{equation}    
\begin{equation} \label{quartic-sigma-omega}     \left\{ \begin{array}{rcl}  \displaystyle  
\sigma_\omega & \, = \,  &  \displaystyle   {{1}\over{12}} \,\, {{N_\omega}\over{D_\omega}} \, \\   
N_\omega & \, = \,  &  \displaystyle 
- 816 \, \sigma_x^4
- 828 \, c_0^2 \, \sigma_x^2  
+ 972 \, c_0^4 \, \sigma_x^2
+ 6480 \, c_0^2 \, \sigma_x^2 \, \sigma_e^2
- 2592  \, c_0^2 \, \sigma_x^4                    \\ &&  \displaystyle   
- 11664 \, c_0^4 \, \sigma_e^2  \, \sigma_x^2
- 816 \, \sigma_e^2 \, \sigma_x^2
+ 1632 \, \sigma_x^3 \, \sigma_e
- 17280 \, \sigma_x^5 \, \sigma_e
+ 13824 \, \sigma_x^4 \, \sigma_e^2           \\ &&  \displaystyle 
+ 3456 \, \sigma_x^6
+ 76 \, \sigma_x^2
+ 51840  \, c_0^2 \, \sigma_x^4  \, \sigma_e^2
- 93312 \, c_0^4 \, \sigma_x^4   \, \sigma_e^2
+ 27 \, c_0^4 - 27 \, c_0^2                           \\ &&  \displaystyle    
+ 4 \, \sigma_e^2
- 324 \, c_0^4 \, \sigma_e^2
+ 180 \, c_0^2 \, \sigma_e^2
- 80 \, \sigma_x \, \sigma_e                        
+ 41472  \, c_0^2 \, \sigma_x^6 
+ 7776 \, c_0^4 \, \sigma_x^4   \\  
D_\omega & \, = \,  &  \displaystyle  \sigma_x \, \big( 
- 192 \, \sigma_x^4
- 828  \, c_0^2 \, \sigma_x^2 
+ 972 \, c_0^4 \, \sigma_x^2
+ 6480 \, c_0^2 \, \sigma_x^2 \, \sigma_e^2
- 2592  \, c_0^2 \, \sigma_x^4                       \\ &&  \displaystyle     
- 11664  \, c_0^4 \, \sigma_e^2 \, \sigma_x^2 
- 192 \, \sigma_e^2 \, \sigma_x^2
+ 384 \, \sigma_x^3 \, \sigma_e
+ 2304 \, \sigma_x^5 \, \sigma_e
+ 4032 \, \sigma_x^4 \, \sigma_e^2                \\ &&  \displaystyle 
- 6336 \, \sigma_x^6
+ 76 \, \sigma_x^2
+ 51840  \, c_0^2 \, \sigma_x^4  \, \sigma_e^2
- 93312  \, c_0^4 \, \sigma_x^4  \, \sigma_e^2
+ 27 \, c_0^4 - 27 \, c_0^2                              \\ &&  \displaystyle 
+ 4 \, \sigma_e^2
- 324 \, c_0^4 \, \sigma_e^2
+ 180 \, c_0^2 \, \sigma_e^2
- 80 \, \sigma_x \, \sigma_e
+ 41472  \, c_0^2  \, \sigma_x^6                     
+ 7776 \, c_0^4 \, \sigma_x^4   \, \big) 
\end{array} \right.    \end{equation}   
For the numerical experiments described above, we have used the following 
values of the equilibrium parameters. 
We did  the  previous theoretical numerical computations 
(\ref{quartic-beta}) 
(\ref{quartic-c2})  
(\ref{quartic-sigma-phi}) 
(\ref{quartic-sigma-beta}) 
(\ref{quartic-sigma-gamma}) 
(\ref{quartic-sigma-chi}) 
(\ref{quartic-sigma-tau})  
(\ref{quartic-sigma-omega}) 
with a 50  %%  40 en 2009 
Digits option in our symbolic manipulation software:
\begin{equation} \label{quartic-equil-numeric}    \left\{ \begin{array}{rcl}  \displaystyle 
%                                                                     valeurs octobre 2009
%  \beta  & \, = \,  &  \displaystyle   2.585949088779921203429311969242201181058 \\  \displaystyle 
% c_0   & \, = \, &   0.5773502691896257645091487805019574556476  \\  \displaystyle  
% c_2  & \, = \, &  0.661536140033197525275388561943564207033   \\  \displaystyle  
%                                                                 fin valeurs octobre 2009
%                                                                     valeurs mai 2010 
%   \beta  & \, = \,  &  \displaystyle   1.649705882352941176470588235294117647059  \\  \displaystyle 
%  c_0   & \, = \, &   0.5773502691896257645091487805019574556476  \\  \displaystyle  
%  c_2  & \, = \, &  1.84375    \\  \displaystyle  
%  \xi & \, = \, &  1
%                                                                 fin valeurs mai 2010 
%                                                         valeurs 26 mai - 7 juin 2010 
s_e            & \, = \, & 0.95057034220532319391634980988593155893536121673004  \\  \displaystyle 
s_x            & \, = \, & 1.8552875695732838589981447124304267161410018552876  \\ \displaystyle   
c_0   & \, = \, &   0.623538  \\  \displaystyle  
c_2  & \, = \, &    2.436118000    \\  \displaystyle  
\xi & \, = \, &  1
%                                                     fin valeurs 26 mai - 7 juin 2010 
\end{array} \right.    \end{equation}   
and the ``quartic'' relaxation parameters: 
\begin{equation} \label{quartic-relax-numeric}    \left\{ \begin{array}{rcl}  \displaystyle 
 \beta  & \, = \,  &  \displaystyle   0.50345521670787922851706691010021052631578947368420  
 \\  \displaystyle
s_\phi         & \, = \, & 0.37925445705024311183144246353322528363047001620745  \\  \displaystyle  
s_\psi         & \, = \, & 1.3 \,\,\,                                            \\  \displaystyle  
s_\varepsilon  & \, = \, & 0.34253657030513141274711235609982461596733955718034  \\  \displaystyle   
s_\xi          & \, = \, & 1.2       \,\,\,                                      \\  \displaystyle  
s_\gamma       & \, = \, & 1.9945477114942149093456286590460711496091258797386   \\  \displaystyle  
s_\chi         & \, = \, & 1.2940799466197218037166471960307116873345869400515   \\  \displaystyle  
s_\tau         & \, = \, & 1.9451616927239606019667013153498030552793202566039   \\  \displaystyle  
s_\omega       & \, = \, & 0.25131560984615404581501005329813711811837082814760   \, . 
%                                                     fin valeurs 26 mai - 7 juin 2010 
\end{array} \right.    \end{equation}   
Of course, we used the previous numbers with the usual truncation allowed by the 32 bit
or 64 bit  floating point arithmetic.   %           computer architectures. 
We observe that for these particular parameters, 
due to the link (\ref{mu-d3q27}) between the shear viscosity $\, \mu \,$ and the 
parameter $\, \sigma_x ,\,$  to   relation  (\ref{zeta-d3q27}) between 
the bulk viscosity $ \, \zeta \,$  and the parameter $\, \sigma_e \,$ 
and to   formula  (\ref{3.3}) making explicit the attenuation  $\, \gamma \,$ 
  of sound waves,  we have 
\begin{equation} \label{quartic-muetal}    \left\{ \begin{array}{rcl}  \displaystyle 
%                                                                     valeurs juin 2010 
\mu            & \, = \, & 0.013    \,\, \lambda \, \Delta x  \\  \displaystyle 
\zeta          & \, = \, & 0.0920492  \,\, \lambda \, \Delta x  \\  \displaystyle 
\gamma         & \, = \, & 0.0546913  \,\, \lambda \, \Delta x   \, . 
%                                                               fin valeurs 12 juin 2010 
\end{array} \right.    \end{equation}    

\bigskip \noindent  $\bullet$ \quad   
The results of simple periodic waves are displayed in Figures 7, 8 and 9.
%  \ref{f-d3q27-stokes}, \ref{f-d3q27-vson}  \ref{f-d3q27-disson}. 
The numerical results confirm that 
parameters proposed  in (\ref{quartic-equil-numeric}) and (\ref{quartic-relax-numeric})
allow one to get fourth order (relative) accuracy.      %%%%G
%

%%%%%%%%%%%%%%%%%%%%%%%%%%%%%%%%%     Figure 07      %%%%%%%%%%%%%%%%%%%%%%%%%%%%%% 
\bigskip     
\centerline { \includegraphics[width=.65  \textwidth]{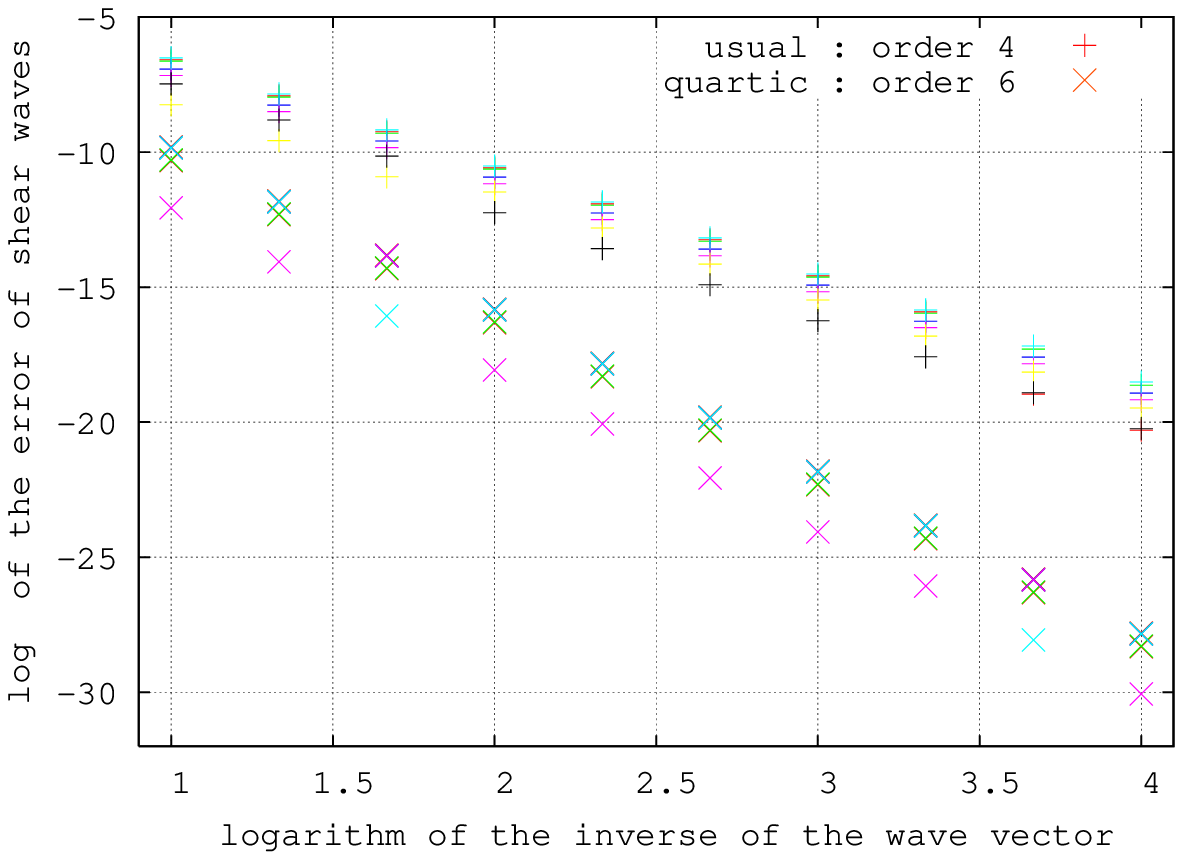}   }
\smallskip \hangindent=7mm \hangafter=1 \noindent {\bf Figure 7}.   \quad  {\it      
Periodic wave with the lattice Boltzmann scheme D3Q27. 
Error for shear eigenvalue $ \,\Gamma_t \,$ 
defined in  (\ref{ondecisail}). For quartic parameters, we have 
$  \displaystyle \, \Gamma_t  -  \nu \, \mid \!\! {\bf k} \!\! \mid^2 =
 {\rm O}(\mid \!\! {\bf k} \!\! \mid^6) . $  }
% \label{f-d3q27-stokes}  \end{figure}  
\bigskip    
%%%%%%%%%%%%%%%%%%%%%%%%%%%%%%%%%%%%%%%%%%%%%%%%%%%%%%%%%%%%%%%% 
%

%%%%%%%%%%%%%%%%%%%%%%%%%%%%%%%%%     Figure 08      %%%%%%%%%%%%%%%%%%%%%%%%%%%%%% 
\bigskip    
\centerline { \includegraphics[width=.65  \textwidth]{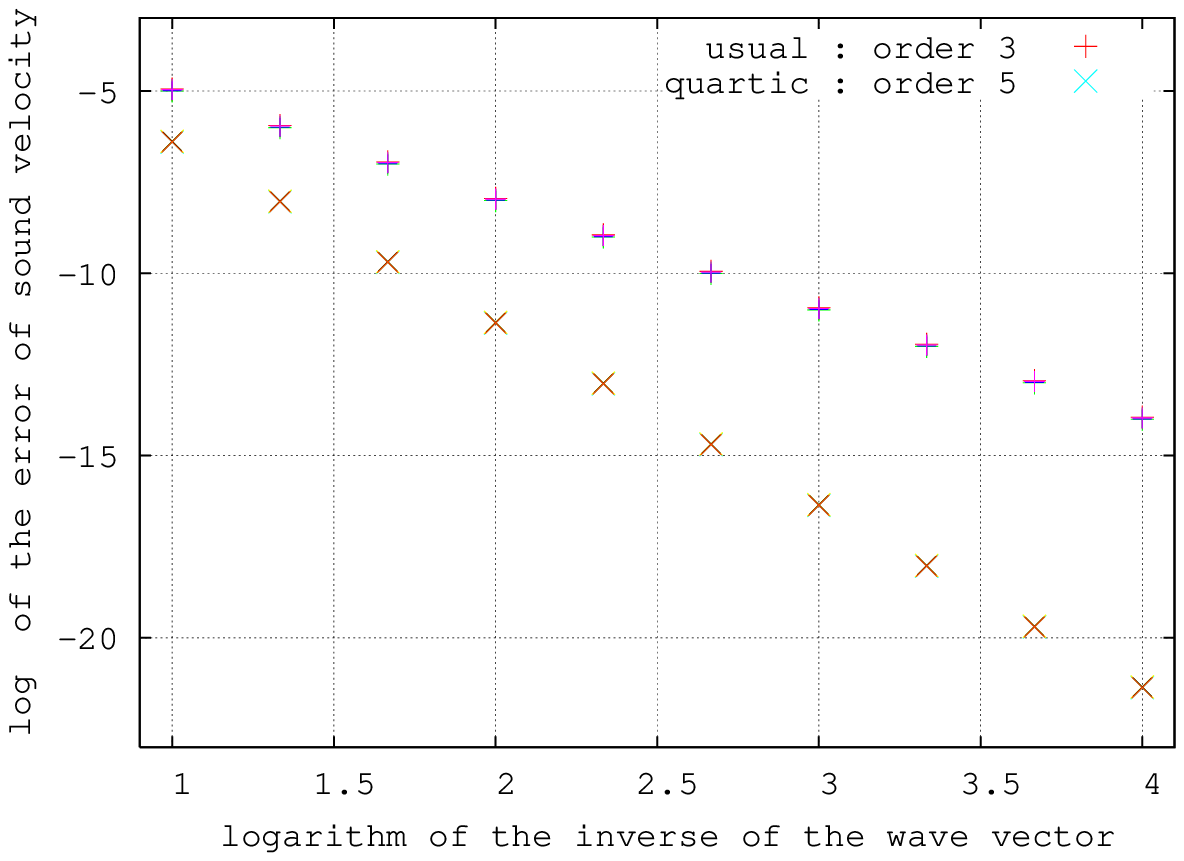}   } 
\smallskip \hangindent=7mm \hangafter=1 \noindent {\bf Figure 8}.   \quad   {\it    
 Periodic wave with the lattice Boltzmann scheme D3Q27. 
Error for the imaginary part of the 
acoustic eigenvalue $ \,\Gamma_\ell  \,$ defined in  (\ref{ondacous}). 
 For quartic parameters, we have  fifth order accuracy: 
 $  \displaystyle \,  {\rm Im} \big(  \Gamma_\ell \big) - c_0 \, \mid \!\! {\bf k} \!\! \mid 
\, \big( 1 - \, \gamma^2 \,  \mid \!\! {\bf k} \!\! \mid^2   / ( 2 \, c_0^2) \big) 
\, = \, {\rm  O}    (  \mid \!\! {\bf k} \!\! \mid^5 )  \,  $ 
 with $ \, \gamma \,$ detailed at relation  (\ref{3.3}).  }
% \label{f-d3q27-vson}  \end{figure}  
\bigskip  
%%%%%%%%%%%%%%%%%%%%%%%%%%%%%%%%%%%%%%%%%%%%%%%%%%%%%%%%%%%%%%%% 
% 

%%%%%%%%%%%%%%%%%%%%%%%%%%%%%%%%%     Figure 09      %%%%%%%%%%%%%%%%%%%%%%%%%%%%%%  
\bigskip  
\centerline { \includegraphics[width=.65  \textwidth]{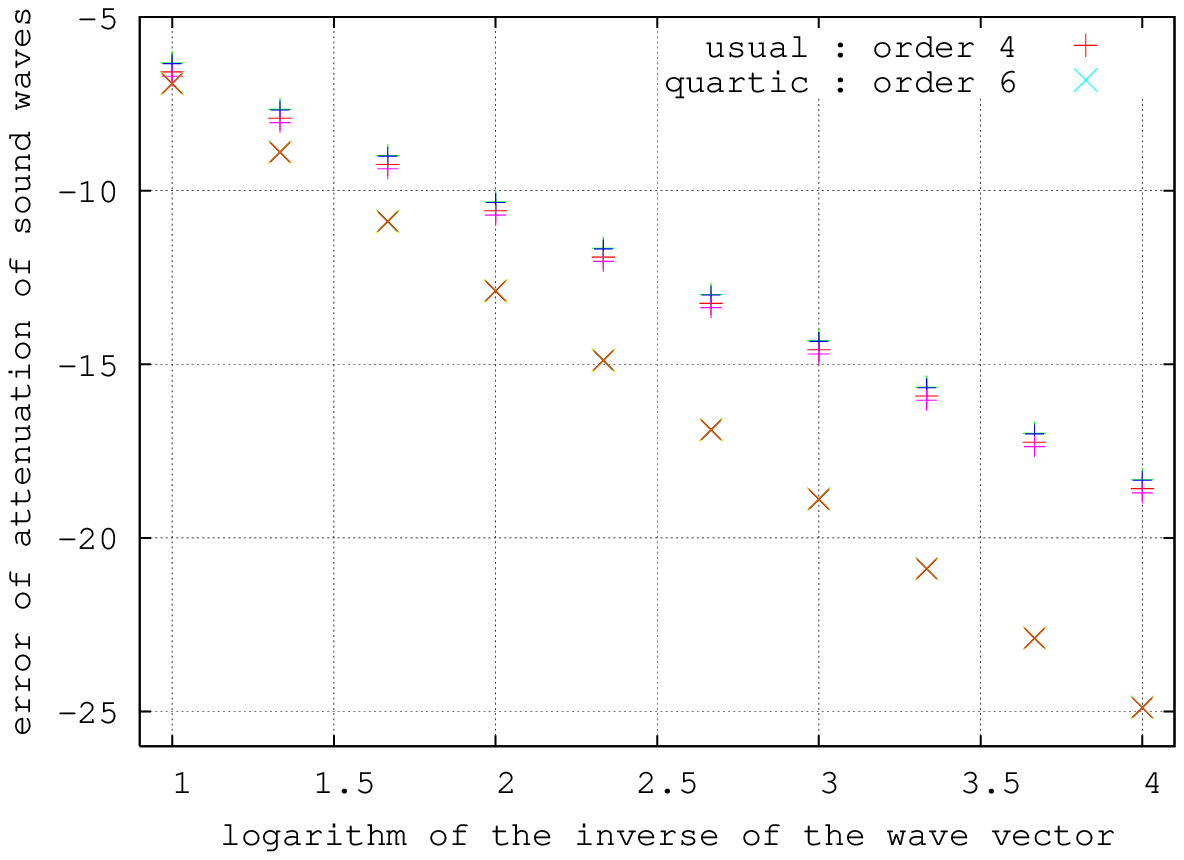}   } 
\smallskip \hangindent=7mm \hangafter=1 \noindent {\bf Figure 9}.   \quad  {\it  
 Periodic wave with the lattice Boltzmann scheme D3Q27. 
Error for the real  part of  acoustic eigenvalue $ \,\Gamma_\ell  . \,$ 
 For quartic parameters, we have 
 $ \displaystyle \,  {\rm Re} \big(  \Gamma_\ell \big) - \gamma \, 
% \Big( {{\zeta}\over{2}} + {{2 \, \mu}\over{3}} \Big) \, 
\mid \!\! {\bf k} \!\! \mid^2   
\, = \,  {\rm O}   (  \mid \!\! {\bf k} \!\! \mid^6 ) \, .  $  } 
% \label{f-d3q27-disson}  \end{figure}  
\bigskip  \bigskip  
%%%%%%%%%%%%%%%%%%%%%%%%%%%%%%%%%%%%%%%%%%%%%%%%%%%%%%%%%%%%%%%% 
%

The numerical experiments have been done for a three-dimensional converging acoustic
wave. 
A pulsating sphere of radius $ \, R = 46.08 \,$ is embedded in a $ \, 95^3 \, $ cube.  
The simulations used only 64  bit or even
32 bit data with CUDA: (using a Nvidia 9800GT card, one full update of D3Q27
takes approximately 11 nanoseconds per site).  
The ``linear anti-bounce-back'' 
numerical boundary conditions (in the sense given in Bouzidi {\it et al.} \cite{bfl}) 
are used to impose a boundary density on the sphere oscillating with a period $ \, T=10 $. 
Since the sound velocity is around $\, 0.62  \,$ (see (\ref{quartic-equil-numeric})), there
are around six mesh points by wave length. In Figure~10,  % \ref{p-d3q27-sphere}, 
the scatter plot  of density {\it   vs} radius after $\, 82 \,$ iterations 
is displayed. The results are of excellent quality with the use of quartic
parameters. Nevertheless, since the attenuation of the sound wave is relatively important 
for  the above parameters 
($ \gamma \approx   0.0547    $ as recalled in  (\ref{quartic-muetal})), 
  the range of propagation is relatively limited (five wave
lengths in our case, see   Figure 10). %   \ref{p-d3q27-sphere}). 
For the less stringent ``isotropic'' case, {\it i.e.} to fix the ideas 
\begin{equation} \label{trt-1}    \left\{ \begin{array}{rcl}  \displaystyle 
\sigma_\phi &=& \displaystyle  {{1}\over{6 \, \sigma_x}}  \\  
\beta  &=& 4 - 9 c_0^2    \\ 
c_0  &=&   \displaystyle  {{1}\over{\sqrt{3}}}
\end{array} \right.    \end{equation} 
we can find situations  of ``TRT type'' satisfying
\begin{equation} \label{trt-2}    \left\{ \begin{array}{rcl}  \displaystyle 
\sigma_\phi &=&  \sigma_\psi   \,\,=\,\,     \sigma_\tau   \,\,=\,\,   \sigma_\omega   \\  
\sigma_e &=&     \sigma_x  \,\,=\,\,  \sigma_\epsilon  \,\,=\,\,  \sigma_\xi  \,\,=\,\,
\sigma_\gamma  \,\,=\,\, \sigma_\chi   
\end{array} \right.    \end{equation}    
with much lower attenuation of sound.
An example is given in figure 10 %  \ref{p-d3q27-sphere}         %   (ddisad9.gz) 
 with parameters of the scheme following  relations (\ref{trt-1}) and  (\ref{trt-2})
with $ \, \sigma_x = 0.006 $, that is 
shear viscosity 0.002    and sound attenuation 0.002. 
These results show that  choosing the goal in terms of accuracy or small attenuation
 will influence the choice of the parameters in the lattice Boltzmann  simulation.

%%%%%%%%%%%%%%%%%%%%%%%%%%%%%%%%%     Figure 10      %%%%%%%%%%%%%%%%%%%%%%%%%%%%%% 
\bigskip  
  \centerline { \includegraphics[width=.65  \textwidth]{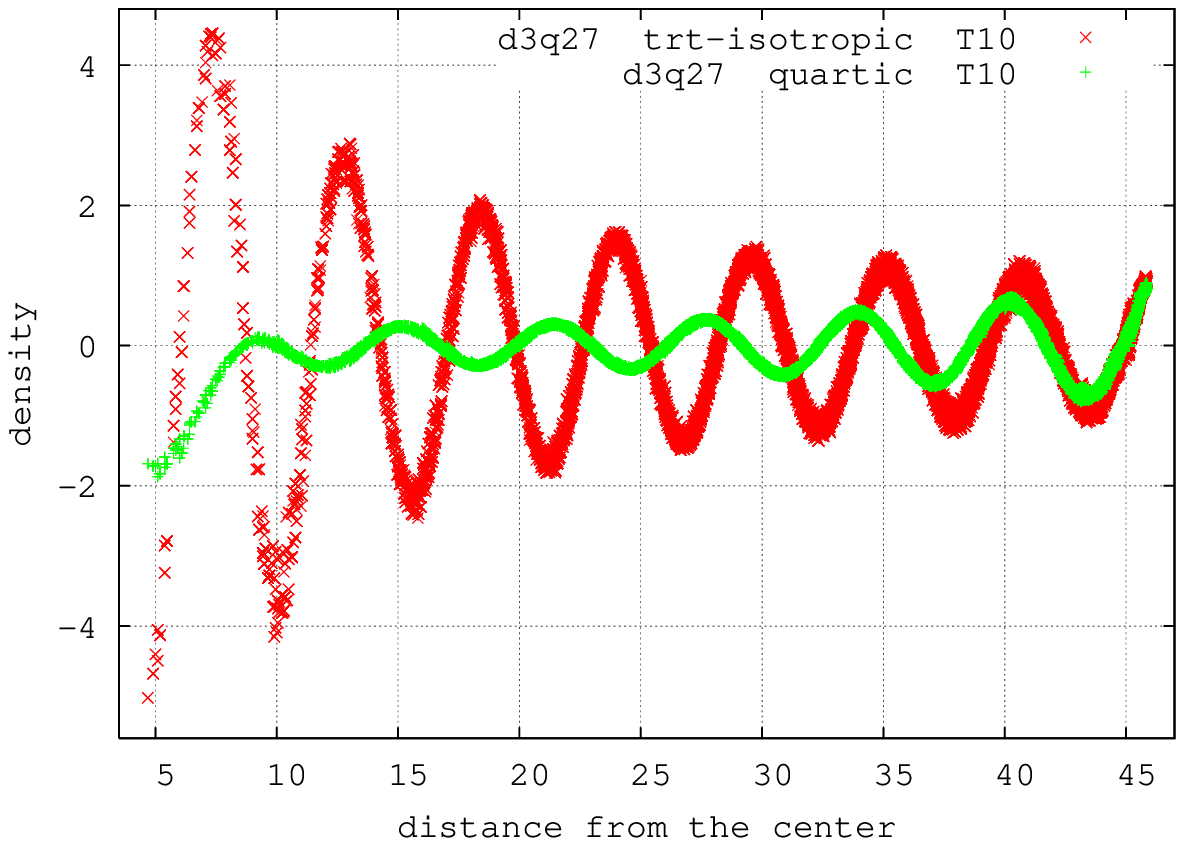}   } 
\smallskip \hangindent=7mm \hangafter=1 \noindent {\bf Figure 10}.   \quad {\it 
 Sound wave emitted by a sphere. Numerical  results with the D3Q27 lattice 
Boltzmann scheme with ``TRT'' isotropic parameters described at relations (\ref{trt-1})
(\ref{trt-2}) and with quartic parameters given by relations (\ref{quartic-equil-numeric})
and  (\ref{quartic-relax-numeric}). Six points by wave length. }  % The density of vertices
% is an increasing function of the radius because all the mesh points are presented. }     
%   \label{p-d3q27-sphere}  \end{figure}  
\bigskip 
%%%%%%%%%%%%%%%%%%%%%%%%%%%%%%%%%%%%%%%%%%%%%%%%%%%%%%%%%%%%%%%% 

%%%%%%%%%%%%%%%%%%%%%%%%%%%%%%%%%%%%%%%%%%%%%%%%%%%%%%%%%%%%%%%%%%%%%%%%%%%%%%%%%%%%%
%%%%%%%%%%%%%%%%%%%%%%%%%%%%%%%%%%%%%%%%%%%%%%%%%%%%%%%%%%%%%%%%%%%%%%%%%%%%%%%%%%%%%
\bigskip \bigskip  \noindent {\bf \large 7) \quad Conclusion}   
%%%%%%%%%%%%%%%%%%%%%%%%%%%%%%%%%%%%%%%%%%%%%%%%%%%%%%%%%%%%%%%%%%%%%%%%%%%%%%%%%%%%%
%%%%%%%%%%%%%%%%%%%%%%%%%%%%%%%%%%%%%%%%%%%%%%%%%%%%%%%%%%%%%%%%%%%%%%%%%%%%%%%%%%%%%
 
\noindent   $\bullet$ \quad  
We have proposed to use the Taylor expansion method in conjunction with the 
use of symbolic manipulation to increase the accuracy of the lattice Boltzmann scheme in the 
multiple time relaxation approach proposed by D.~d'Humi\`eres  for linear acoustic waves. 
The result is the use of the previous scheme with very particular ``quartic'' parameters.
We have obtained  a family of such parameters for the D3Q27 numerical scheme. 
Numerical experiments confirm  the predictions of the theoretical analysis.  
The catch is to make sure that situations are stable. This problem cannot be solved in    %%%%G
practice with developments around small wave numbers. 
%%%%%%%%%%%%%%%%%%%%%%%%%%%%%%%%%%%%%%%%%%%%%%   modif Pierre 14 juin 2010 
%% We have searched at random in the  %%%%G
%% parameter space to the get the values shown above.
%% More efficient techniques to analyze
%% the stability of the lattice Boltzmann scheme are needed. 
%%%%%%%%%%%%%%%%%%%%%%%%%%%%%%%%%%%%%%%%%%%%%%   fin modif Pierre 14 juin 2010  
In the absence of efficient techniques to find stable situations with small viscosities 
and/or sound attenuation, we tried at random a number of sets of the free parameters and have
used the best of them for the explicit results shown above. 
We note that the less stringent ``isotropic''  condition  is
compatible with small viscosities which may be quite useful for practical
simulations. Comparison of various stencils is under way and will be
presented elsewhere.

%%%%%%%%%%%%%%%%%%%%%%%%%%%%%%%%%%%%%%%%%%%%%%%%%%%%%%%%%%%%%%%%%%%%%%%%%%%%%%%%%%%%%
%%%%%%%%%%%%%%%%%%%%%%%%%%%%%%%%%%%%%%%%%%%%%%%%%%%%%%%%%%%%%%%%%%%%%%%%%%%%%%%%%%%%%
\bigskip \bigskip   \noindent {\bf \large Acknowledgments}   
%%%%%%%%%%%%%%%%%%%%%%%%%%%%%%%%%%%%%%%%%%%%%%%%%%%%%%%%%%%%%%%%%%%%%%%%%%%%%%%%%%%%%
%%%%%%%%%%%%%%%%%%%%%%%%%%%%%%%%%%%%%%%%%%%%%%%%%%%%%%%%%%%%%%%%%%%%%%%%%%%%%%%%%%%%%

\noindent 
The referees conveyed  to the authors very interesting  remarks.
Some of them  have been incorporated into  the present edition of the article.

%%%%%%%%%%%%%%%%%%%%%%%%%%%%%%%%%%%%%%%%%%%%%%%%%%%%%%%%%%%%%%%%%%%%%%%%%%%%%%%%%%%%%
%%%%%%%%%%%%%%%%%%%%%%%%%%%%%%%%%%%%%%%%%%%%%%%%%%%%%%%%%%%%%%%%%%%%%%%%%%%%%%%%%%%%%
\bigskip \bigskip  \noindent {\bf \large  References } 
%%%%%%%%%%%%%%%%%%%%%%%%%%%%%%%%%%%%%%%%%%%%%%%%%%%%%%%%%%%%%%%%%%%%%%%%%%%%%%%%%%%%%
%%%%%%%%%%%%%%%%%%%%%%%%%%%%%%%%%%%%%%%%%%%%%%%%%%%%%%%%%%%%%%%%%%%%%%%%%%%%%%%%%%%%%

\medskip

\end{document}